\documentclass[a4paper,reqno,12pt]{amsart}

\usepackage{bbm}

\usepackage{wrapfig}
\usepackage{tikz}
\usetikzlibrary{calc}
\usepackage{pgfplots}
\usepackage{ucs}
\usepackage[utf8x]{inputenc}
\usepackage[T1]{fontenc}

\usepackage{amsmath}
\usepackage{amsfonts}
\usepackage{amssymb}
\usepackage{amsthm}
\usepackage{xcolor}

\usepackage{paralist}
\usepackage{enumitem}

\usepackage[margin=1cm,labelfont=sc,font={small}]{caption}

\usepackage{graphicx}
\DeclareGraphicsExtensions{.png,.pdf,.jpg}
\graphicspath{{pics/}}

\theoremstyle{plain}
\newtheorem{theorem}{Theorem}[section]

\newtheorem{conjecture}[theorem]{Conjecture}
\newtheorem{definition}{Definition}

\theoremstyle{definition}

\newcommand{\mintwo}[2]{\min_{\substack{#1 \\ #2}}} 
\newcommand{\sumtwo}[2]{\sum_{\substack{#1 \\ #2}}} 
\newcommand{\abs}[1]{\left| #1\right|}

\newcommand{\calA}{\mathcal{A}}
\newcommand{\calB}{\mathcal{B}}
\newcommand{\calC}{\mathcal{C}}

\newcommand{\calH}{\mathcal{H}}

\newcommand{\calL}{\mathcal{L}}

\newcommand{\calP}{\mathcal{P}}

\newcommand{\calS}{\mathcal{S}}

\newcommand{\calW}{\mathcal{W}}

\newcommand{\calY}{\mathcal{Y}}

\newcommand{\fra}{\mathfrak{a}}
\newcommand{\frb}{\mathfrak{b}}

\newcommand{\frn}{\mathfrak{n}}

\newcommand{\bbA}{\mathbb{A}}
\newcommand{\bbB}{\mathbb{B}}

\newcommand{\bbE}{\mathbb{E}}

\newcommand{\bbH}{\mathbb{H}}

\newcommand{\bbL}{\mathbb{L}}

\newcommand{\bbN}{\mathbb{N}}

\newcommand{\bbP}{\mathbb{P}}
\newcommand{\bbQ}{\mathbb{Q}}
\newcommand{\bbR}{\mathbb{R}}
\newcommand{\bbS}{\mathbb{S}}

\newcommand{\bbX}{\mathbb{X}}

\newcommand{\bbZ}{\mathbb{Z}}

\newcommand{\sfa}{\mathsf a}
\newcommand{\sfb}{\mathsf b}

\newcommand{\sfe}{\mathsf e}
\newcommand{\sff}{\mathsf f}
\newcommand{\sfg}{\mathsf g}
\newcommand{\sfh}{\mathsf h}

\newcommand{\sfl}{{\mathsf l}}

\newcommand{\sfn}{{\mathsf n}}

\newcommand{\sfp}{{\mathsf p}}

\newcommand{\sfr}{{\mathsf r}}

\newcommand{\sfv}{{\mathsf v}}
\newcommand{\sfw}{{\mathsf w}}
\newcommand{\sfx}{{\mathsf x}}
\newcommand{\sfy}{{\mathsf y}}
\newcommand{\sfz}{{\mathsf z}}

\newcommand{\sfA}{\mathsf{A}}
\newcommand{\sfB}{\mathsf{B}}
\newcommand{\sfC}{\mathsf{C}}
\newcommand{\sfD}{\mathsf{D}}

\newcommand{\sfG}{\mathsf{G}}
\newcommand{\sfH}{\mathsf{H}}

\newcommand{\sfK}{\mathsf{K}}
\newcommand{\sfL}{\mathsf{L}}

\newcommand{\sfP}{\mathsf{P}}
\newcommand{\sfQ}{\mathsf{Q}}

\newcommand{\sfS}{\mathsf{S}}
\newcommand{\sfT}{\mathsf{T}}

\newcommand{\sfW}{\mathsf{W}}
\newcommand{\sfX}{\mathsf{X}}
\newcommand{\sfY}{\mathsf{Y}}
\newcommand{\sfZ}{\mathsf{Z}}

\newcommand{\ul}[1]{\underline{#1}}

\newcommand{\ur}{\underline{r}}

\newcommand{\uv}{\underline{v}}
\newcommand{\uu}{\underline{u}}

\newcommand{\ugamma}{\underline{\gamma}}

\newcommand{\ueta}{\underline{\eta}}

\newcommand{\uphi}{\underline{\varphi}}

\newcommand{\uB}{\underline{B}}
\newcommand{\uR}{\underline{R}}
\newcommand{\uX}{\underline{X}}

\newcommand{\Rd}{\bbR^d}

\newcommand{\Zd}{\bbZ^d}

\newcommand{\setof}[2]{\left\{#1 \,:\, #2 \right\}}

\newcommand{\defby}{\stackrel{\text{\tiny{\rm def}}}{=}}

\newcommand{\IF}[1]{\mathbf{1}_{\{#1\}}}

\newcommand{\lb}{\left(}
\newcommand{\rb}{\right)}
\newcommand{\lbr}{\left\{}
\newcommand{\rbr}{\right\}}

\newcommand{\dd}{{\rm d}}

\newcommand{\step}[1]{S{\small TEP}\,#1.}

\newcommand{\1}{\mathbbm{1}}

\newcommand{\smo}[1]{{\mathrm o}\lb #1\rb }

\newcommand{\df}{\stackrel{\Delta}{=}}

\newcommand{\be}[1]{\begin{equation}\label{#1}}
\newcommand{\ee}{\end{equation}}

\newcommand{\Hla}{H_\lambda}

\newcommand{\bsetof}[2]{\bigl\{#1\,:\,#2\bigr\}}

\newcommand{\normII}[1]{\|#1\|_{\scriptscriptstyle 2}}
\newcommand{\normsup}[1]{\|#1\|_{\scriptscriptstyle\infty}}

\newcommand{\uvec}{\boldsymbol{e}}

\newcommand{\Ham}{\calH}	
\newcommand{\PHS}{\bbH_+}	
\newcommand{\NHS}{\bbH_-}	
\renewcommand{\emptyset}{\varnothing}
\newcommand{\betac}{\beta_{\mathrm{c}}}
\newcommand{\hw}{h_{\mathrm{w}}}
\newcommand{\bfn}{\boldsymbol{v}}
\newcommand{\eig}{\sfe}

\author{Dmitry Ioffe}
\address{Faculty of IE\&M, Technion, Haifa 32000, Israel}
\email{ieioffe@ie.technion.ac.il}
\thanks{DI was supported by the Israeli Science Foundation grant 
1723/14 and  by The Leverhulme Trust through International
Network Grant {\em Laplacians, Random Walks, Bose Gas, Quantum Spin Systems}.}
\author{Yvan Velenik}
\address{Section de Math\'ematiques, Universit\'e de Gen\`eve, 
CH-1211 Gen\`eve, Switzerland}
\email{yvan.velenik@unige.ch}
\thanks{YV was partially supported by the Swiss National Science Foundation and the NCCR SwissMAP.}

\title[Low-temperature interfaces and Ferrari--Spohn diffusions]
{Low-temperature interfaces: Prewetting, layering, faceting  and Ferrari--Spohn 
diffusions}

\begin{document}
\maketitle
\begin{abstract}
In this paper, we survey and discuss various surface phenomena such 
as prewetting, layering and faceting for a family of two- and three-dimensional
low-temperature models of statistical mechanics, notably Ising models
 and $2+1$-dimensional solid-on-solid (SOS) models,  with a particular accent on 
 scaling regimes which lead or, in most cases, are conjectured to
 lead to Ferrari--Spohn type diffusions.  
\end{abstract}

\section{Introduction and structure of the paper} 

Ferrari--Spohn (FS) diffusions were introduced in~\cite{FS05} as 
cube-root scaling limits of Brownian motion constrained
to stay above circular and parabolic barriers. It turns out that  
FS diffusions and their Dyson counterparts are the universal scaling 
limits for a class of random walks~\cite{IoShVe2015} or, respectively, 
ordered random walks~\cite{IVW16} under generalized area tilts. 

Random walks under vanishing area tilts naturally appear as effective
descriptions for either phase separation lines in two-dimensional models
of Statistical mechanics in regimes related to  critical prewetting, 
or for level lines of random low-temperature surfaces in three dimensions,
in the regime when the number of layers grows and/or when the linear
size of the system goes to infinity. 

Making the above statement precise sets up the stage for a host of open
problems, some of them --- for instance  those which arise in the context 
of Ising interfaces in three dimensions --- seem, for the moment, notoriously 
difficult. On the other hand, phase separation lines in two dimensions
and level lines of SOS interfaces are more tractable objects, and we expect
progress in understanding FS scaling for the latter in a foreseeable future.

The paper is organized as follows: In Section~\ref{sec:wetting} we try 
to survey the existing rigorous results and many remaining unsolved issues
on wetting transition, mostly in context of Ising models on \(\Zd\). 
Critical prewetting, which is the somehow central notion of the 
scaling phenomena in question, is discussed in Section~\ref{sec:crit-pre}.
Various types of 2+1-dimensional SOS-type models are discussed in 
Section~\ref{sec:Effective}, including models with hard-wall constraint
(entropic repulsion), bulk and boundary fields and models coupled with 
high- and low-density Bernoulli fields below and above the interfaces. 
 
Rigorous results on scaling of non-colliding random walks under area 
tilts are collected and briefly explained in Section~\ref{sec:RW}. 
 
Definitions and properties of Wulff shapes, the definition of FS and
Dyson FS diffusions, as well as a sketch of a strategy to prove FS-scaling
for self-avoiding walks under area tilts are relegated to
Appendix~\ref{app:Wulff}, Appendix~\ref{app:FS} and Appendix~\ref{app:SAW}
respectively.

\section{The wetting transition}
\label{sec:wetting}

Consider a system with (at least) two thermodynamically stable phases $A$ and 
$B$. We assume that the bulk of the system is occupied by phase $A$ and that 
the latter interacts with a substrate favoring phase $B$. In such a 
situation, depending on how much preference the substrate displays for phase 
$B$, two different scenarios are possible:

\medskip
\begin{minipage}{.75\textwidth}
\noindent
$\bullet$ The bulk phase is in contact with the substrate, with only microscopic 
droplets of phase $B$ attached to the latter. This is called the regime of 
\emph{partial wetting} (of the substrate by phase $B$).
\end{minipage}
\begin{minipage}[r]{3cm}
\centering
\includegraphics[width=2.5cm]{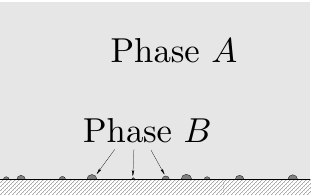}
\end{minipage}

\begin{minipage}{.75\textwidth}
\noindent
$\bullet$
A mesoscopic layer of phase $B$ covers the substrate, thus separating 
the latter from the bulk phase. This is known as the regime of \emph{complete wetting} (of the substrate by phase $B$).
\end{minipage}
\begin{minipage}[r]{3cm}
\centering
\includegraphics[width=2.5cm]{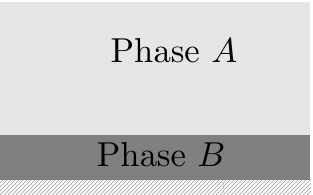}
\end{minipage}

\medskip

The phase transition from one behavior to the other, usually triggered by a 
change of temperature (although it is often more convenient mathematically to 
use a change in the interaction between the substrate and the phases) is 
known as the \emph{wetting transition}; the latter is a prime example of a  
\emph{surface phase transition}.

\medskip
In the remainder of this section, we discuss the phenomenon of wetting in 
some specific two- and three-dimensional systems.
In order to keep the discussion as concrete as possible, we limit most of it 
to the Ising model and, in Section~\ref{sec:Effective}, 
to the corresponding effective interface models.
We shall also make a few remarks on the analogous phenomena occurring in the 
Blume--Capel model. The general phenomenology described here has, of course, 
a much broader domain of validity.

\subsection{The wetting transition in the Ising model}

\begin{figure}[t]
\begin{center}
\includegraphics{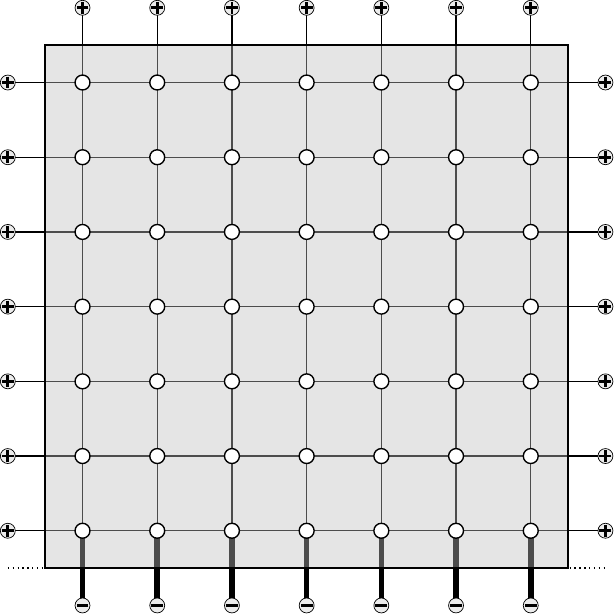}
\end{center}
\caption{The box $\Delta_3$ with $\pm$ boundary condition (in \(\bbZ^2\)). 
The thick edges indicate the modified interaction between the spins in the 
bottom row of $\Delta_3$ and their neighbor in the substrate.}
\label{fig:IsingBox}
\end{figure}

Let 
\[
\PHS^d\defby\setof{i=(i_1,\ldots,i_d)\in\Zd}{i_d\geq 0} 
\]
and 
$\NHS^d\defby\Zd\setminus\PHS^d$. Let $\beta,h\geq 0$ and set
\be{eq:J-h}
J_{i,j} \defby
\begin{cases}
1		& \text{if } \{i,j\}\subset\PHS^d,\\
h		& \text{otherwise.}
\end{cases}
\ee
We consider the nearest-neighbor ferromagnetic Ising model in the box 
\[
\Delta_n \defby \{-n,\ldots,n\}^{d-1}\times\{0,\ldots 2n\}
\]
with Hamiltonian (see Figure~\ref{fig:IsingBox})
\begin{equation}
\label{eq:Ham:Ising2dNoBulkField}
\Ham_{n;h}(\sigma)
\defby
- \sum_{\{i,j\}\cap\Delta_n\neq\emptyset}  J_{i,j} \sigma_i\sigma_j.
\end{equation}
In other words, the interaction between the spins located on the bottom 
``face'' of \(\Delta_n\) and their neighbor below them is modified; this 
models the interaction of the system with a substrate.

We consider two sets of configurations:
\be{eq:pm-bc}
\Omega_n^\pm
\defby
\bsetof{\sigma\in\{\pm 1\}^{\Zd}}{\sigma_i=1 \text{ if } 
i\in\PHS^d\setminus\Delta_n \text{ and } \sigma_i=-1 \text{ if } 
i\in\NHS^d}
\ee
and
\[
\Omega_n^-
\defby
\bsetof{\sigma\in\{\pm 1\}^{\Zd}}{\sigma_i=-1 \text{ for all } 
i\not\in\Delta_n}.
\]
The corresponding Gibbs measures are then defined by
\[
\mu^\pm_{n;\beta,h}(\sigma) \defby \IF{\sigma\in\Omega_{\Delta_n}^\pm}
\frac{e^{-\beta\Ham_{n;h}(\sigma)}}{\sfZ_{\Delta_n;\beta,h}^\pm} 
\quad\text{ and }\quad
\mu^-_{n;\beta,h}(\sigma) \defby \IF{\sigma\in\Omega_{\Delta_n}^-}
\frac{e^{-\beta\Ham_{n;h}(\sigma)}}{\sfZ_{\Delta_n;\beta,h}^-} 
.
\]

We shall always assume that $\beta>\betac(d)$, where 
$\betac(d)\in (0,\infty)$ is the inverse critical temperature 
of the model. In this regime, typical configurations (for large $n$) are best 
characterized in terms of the corresponding Peierls contours, that is, the 
maximal connected components of the boundary of the set
\[
\bigcup_{i\in\Zd:\sigma_i=-1} \setof{x\in\Rd}{\normsup{x-i}\leq 1/2}.
\]

Note that, under $\mu^-_{n;\beta,h}$, $-$ spins are favored along the inner boundary
of the box $\Delta_n$, which always result in the $-$ phase filling the box. The
situation is more interesting under $\mu^\pm_{n;\beta,h}$, since then the bottom
side (corresponding to the substrate) favors $-$ spins, while the other sides favor
$+$ spins. It is this competition that will lead to the wetting transition.

\subsubsection{Entropic repulsion when \(h=1\).}

Let us first consider the simpler problem in which \(h=1\). In this case, 
under the measure \(\mu^\pm_{n;\beta,h=1}\), the boundary condition induces an 
interface and it is natural to investigate what effect the bottom wall has on the
interface.
In~\cite{FrPf1987b}, it is proved that, at any temperature and in any dimension, 
the interface is repelled infinitely far away from the wall as \(n\to\infty\), in
the sense that the weak limits of \(\mu^-_{n;\beta,h=1}\) and \(\mu^\pm_{n;\beta,h=1}\) 
as \(n\to\infty\) coincide (which means that the layer of \(-\) 
phase along the wall extends infinitely far away from the wall). 
In the same paper, this was shown to imply also a weak form of 
delocalization of the interface at sufficiently low temperature. 
In dimensions \(d\geq 3\) and at low enough temperatures, 
these delocalization results can be considerably strengthened~\cite{HoZa}, 
proving that ``most'' of the interface lies at a height of order \(\log n\) 
above the wall. This phenomenon is known as \emph{entropic repulsion}.

The issue of wetting becomes now clear: when \(|h| < 1\), it becomes energetically
favorable for the interface to lie along the substrate, and this attraction will
compete with entropic repulsion: 
partial wetting will occur when attraction wins, while complete wetting will occur
when entropic repulsion wins. Let us now discuss these issues in more detail.

\subsubsection{The two-dimensional case.}

In this section, we describe some results that have been obtained concerning 
the wetting transition in the two-dimensional Ising model.

Observe that all configurations in $\Omega_n^\pm$ possess a unique Peierls 
contour $\gamma$ of infinite length. In fact, with 
$\mu^\pm_{n;\beta,h}$-probability going to $1$ as $n\to\infty$, all other 
contours are of diameter at most $K(\beta)\log n$ for some $K(\beta)<\infty$.

The behavior of the contour $\gamma$ depends very strongly on the value 
of the parameter $h$. Namely, there exists $\hw=\hw(\beta)\geq 0$ such that, 
with $\mu^\pm_{n;\beta,h}$-probability going to $1$ as $n\to\infty$, the 
following occurs (see Figure~\ref{fig:wettingIsing2d}):

\begin{figure}[t]
\begin{center}
\includegraphics[width=\textwidth]{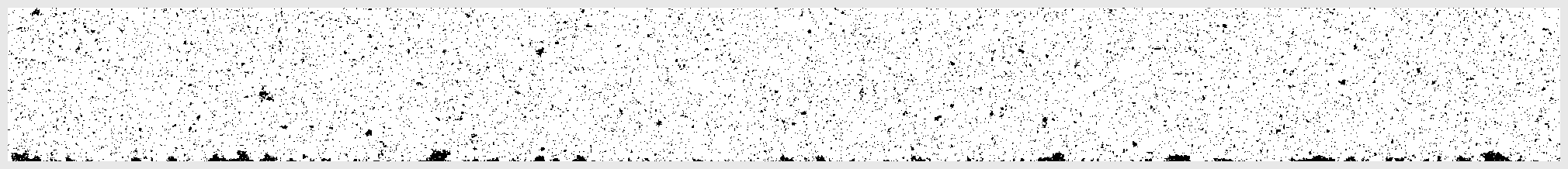}
\centerline{\small Partial wetting (\(0\leq h\leq\hw\))}

\smallskip
\includegraphics[width=\textwidth]{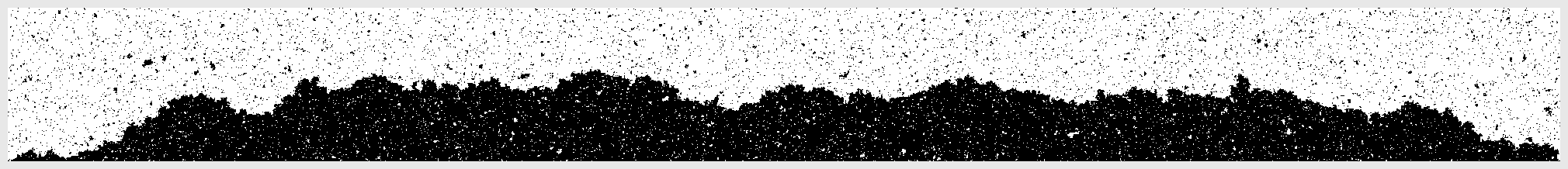}
\centerline{\small Complete wetting (\(h\geq\hw\))}
\end{center}
\caption{Two-dimensional Ising model at $\beta=0.5$; $+$ spins are white, $-$ 
spins black. \textit{Top:} partial wetting regime ($h=0.3$). \textit{Bottom:} 
complete wetting regime ($h=1$).}
\label{fig:wettingIsing2d}
\end{figure}%

\begin{itemize}
\item[] \textbf{Partial wetting:} When $h\in[0,\hw)$, the Hausdorff distance 
between $\gamma$ and the line 
$\calL\defby\setof{x=(x_1,x_2)\in\bbR^2}{x_2=-\tfrac12}$ is of order 
$\log n$. More precisely, the diameter of the maximal connected components of 
$\gamma\setminus\calL$ have exponential tails (see 
Fig.~\ref{fig:diamondsOnTheWall}); this follows from a 
combination of the methods in~\cite{PfVe1997,PfVe99} and~\cite{FrIoVe}.
\begin{figure}[ht]
\begin{center}
\includegraphics[width=13cm]{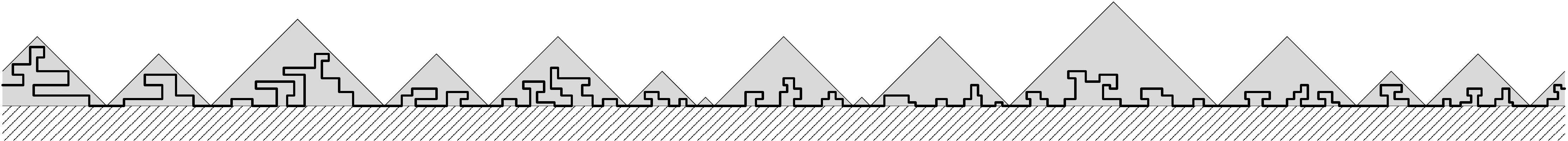}
\end{center}
\caption{Typical realizations of the interface of the two-dimensional Ising 
model in the partial wetting regime can be split into pieces 
included in disjoint triangles with their basis on the wall. The length of 
the basis of these triangles has exponential tail.}
\label{fig:diamondsOnTheWall}
\end{figure}
\item[] \textbf{Complete wetting:} If $h\geq\hw$, the Hausdorff distance 
between $\gamma$ and the line $\calL$ is of order $n^{1/2}$.
More precisely, it is expected that, under diffusive scaling, $\gamma$ 
converges to a Brownian excursion when $h>\hw$; this should follow from a 
combination of the methods in~\cite{PfVe99,CIV03,GI05, CIL10}.
Moreover, on the basis of a similar result for the corresponding effective 
interface model in~\cite{DGZ2005}, it is natural to conjecture that the 
scaling limit is a modulus of Brownian bridge  when $h=\hw$.
\end{itemize}

\begin{figure}[t]
\begin{center}
\includegraphics{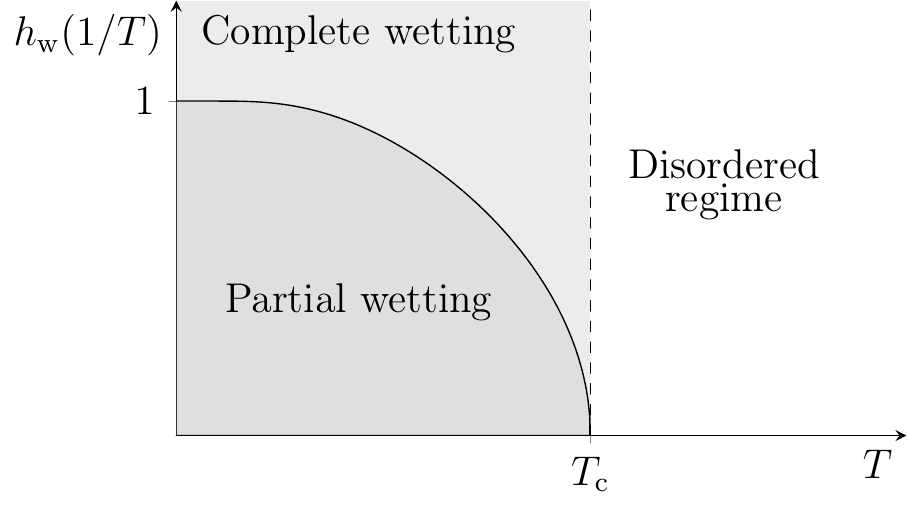}
\end{center}
\caption{The graph of the critical boundary coupling constant $\hw(\beta)$ as 
a function of the temperature $T=1/\beta$, in the two-dimensional Ising 
model. The regime of partial, respectively complete, wetting correspond to 
the dark, respectively light, shaded area. For $T\geq T_{\mathrm{c}}$, there 
is no phase coexistence.}
\label{fig:hw}
\end{figure}%

Remarkably, $\hw$ can be determined explicitly in this case~\cite{Ab1980}: it 
is given by the unique nonnegative solution (see Figure~\ref{fig:hw}) to
\[
\exp(2\beta)\bigl( \cosh(2\beta) - \cosh(2\beta\hw) \bigr) = 
\sinh(2\beta).
\]
The sudden change in the behavior of $\gamma$ as $h$ crosses the value $\hw$ 
is the manifestation of the wetting transition in this setting.

\subsubsection{The three-dimensional case.}

The current understanding of the wetting transition in the three-dimensional 
Ising model is much more rudimentary. This is mainly due to the fact that 
random lines are easier to describe than random surfaces.
But more is true: the qualitative behavior is substantially more 
complicated in the three-dimensional model, as explained below.

\medskip
One first {complication}  when $d=3$ is the conjectured existence of a 
roughening transition. Namely, it is well known~\cite{Dobrushin} that, at 
sufficiently low temperature, imposing boundary conditions that induce the 
presence of a horizontal interface between $+$ and $-$ bulk phases results in 
an interface that is rigid: the latter coincides with a perfect plane, apart 
from rare local deviations. It is conjectured that this interface loses its 
rigid character at a temperature $\beta_{\rm r}<\betac(3)$. This problem 
remains completely open from a mathematical point of view, only the regime of 
very low temperatures being well understood. As a consequence, all the 
results discussed below are valid only at sufficiently large values of 
$\beta$; in particular, $\beta>\beta_{\rm r}$.

\medskip
Let us review what is known rigorously about the wetting transition in this 
model. First, the existence of a value $\hw(\beta)$ separating regimes of 
partial and complete wetting has been proved in~\cite{FrPf1987b}. Moreover, 
explicit lower and upper bounds on $\hw(\beta)$ have been derived in this 
work (in particular, $\hw(\beta)\leq 1$ for all $\beta$), and a weak form of 
localization, resp.\ delocalization, of the corresponding Peierls contour in 
the partial wetting, resp.\ complete wetting, regime was established.

\medskip
It is however expected that the situation is actually more interesting in 
this case: inside the partial wetting regime, it is conjectured that there is 
an \emph{infinite} sequence of first-order phase transitions, known as 
\emph{layering transitions}, at which the microscopic height of the (rigid) 
interface increases by one unit. The precise behavior is not known, but some 
(not all fully rigorous) quantitative information about the location of the 
first few lines has been obtained in~\cite{Ba88,Ba07,AlDuMi10}.

Rigorous results on layering transitions in the context of effective 
SOS models were obtained in~\cite{AlDuMi11}, see a discussion in the end of 
Subsection~\ref{subsub:SOS} below. 

\subsubsection{Some terminology}
\label{sub-terminology}
Before briefly presenting alternative equivalent characterizations of the 
wetting transition, we need some terminology, which will also be useful later.

\paragraph{\bf Surface tension.}
The free energy per unit-area (unit-length when $d=2$) associated to a linear 
interface is measured by the surface tension. In the \(d\)-dimensional Ising 
model, the latter quantity can be defined as follows. Let 
$\Lambda_n=\{-n,\ldots, n\}^d$ and let $\bfn$ be a unit-length vector in 
$\Rd$. Let also $\bbH^d_{\bfn,+} \defby \setof{i\in\Zd}{i\cdot\bfn\geq 0}$ 
and 
$\bbH^d_{\bfn,-} \defby \Zd\setminus\bbH^d_{\bfn,+}$ be the two half-spaces 
delimited by the plane going though the origin and with normal \(\bfn\), and 
set
\[
\Omega^{\bfn}_n \defby \bsetof{\sigma\in\{\pm 1\}^{\Zd}}{\sigma_i = 
\IF{i\in\bbH^d_{\bfn,+}} - \IF{i\in\bbH^d_{\bfn,-}} \text{ for all 
}i\not\in\Lambda_n}
\]
and $\Omega^+_n \defby \bsetof{\sigma\in\{\pm 1\}^{\Zd}}{\sigma_i = 1 
\text{ for all }i\not\in\Lambda_n}$. The surface tension per unit-area in 
the direction normal to $\bfn$ is defined by
\be{eq:st-Ising}
\tau_\beta(\bfn) \defby \lim_{n\to\infty} - 
\frac{1}{\beta|A_n(\bfn)|} \log 
\frac{\sfZ^{\bfn}_{\Lambda_n;\beta}}{\sfZ^{+}_{\Lambda_n;\beta}},
\ee
where $|A_n(\bfn)|$ is the area (more precisely, the $(d-1)$-dimensional 
Hausdorff measure) of the 
set $\setof{x\in [-n,n]^d}{x\cdot\bfn = 0}$ and, for $\star \in\{+,\bfn\}$,
\[
\sfZ^\star_{\Lambda_n;\beta} \defby \sum_{\sigma\in\Omega^\star_n} 
\exp\bigl\{
\beta\sum_{\{i,j\}\cap\Lambda_n\neq\emptyset} \sigma_i\sigma_j
\bigr\}.
\]
The existence of the limit in the definition of \(\tau_\beta(\bfn)\) is 
proved in~\cite{MeMiRu}.
The function $\tau_\beta$ can be extended to $\Rd$ by positive 
homogeneity. Namely, one sets $\tau_\beta(0) \defby 0$ and, for any $0\neq 
x\in\Rd$, $\tau_\beta(x) \defby \tau_\beta(x/\normII{x})\normII{x}$.
This function can be shown to be convex~\cite{MeMiRu}. Moreover, one can 
show~\cite{BLP80,LP81} that $\tau_\beta$ is an order parameter in the sense 
that \(\tau_\beta > 0\) if and only if \(\beta>\betac(d)\).

Equilibrium crystal shapes (Wulff shapes) associated to the surface tension 
$\tau_\beta$ are described in Appendix~\ref{app:Wulff}. 

In the three-dimensional Ising model at sufficiently low temperature 
(that is, at sufficiently large values of \(\beta\)), it is possible to 
prove that the surface tension \(\tau_\beta(\bfn)\) does not behave smoothly 
as a function of \(\bfn\). Namely, parameterizing \(\bfn\) as a function of 
the polar angle \(\theta\) and azimuthal angle \(\varphi\), we can 
write \(\tau_\beta(\theta,\varphi) \equiv \tau_\beta(\bfn(\theta,\varphi))\). 
It is then possible to prove that 
\(\partial\tau_\beta(\theta,\varphi)/\partial\theta\) is 
discontinuous at \(\theta=0\) for all \(\varphi\). This has an 
impact on the geometry of the Wulff shape, which develops facets 
(that is, flat portions of positive two-dimensional Hausdorff measure) 
orthogonally to each lattice direction. These facets are conjectured to 
disappear at the roughening transition. Proofs and additional information on 
these issues can be found in~\cite{Miracle-Sole1995}. Note that it is still an 
open problem to determine whether these \emph{macroscopically} flat pieces of 
the Wulff shape  are actually also \emph{microscopically} flat 
(in the same sense as the Dobrushin interface). 
We will see later, in Section~\ref{sec:Effective}, what can be said about that 
in the simpler context of effective interface models.

\paragraph{\bf Wall free energy.}
We now want to introduce a quantity measuring the free energy per unit-area
associated to an interface along a substrate. The definition is very similar 
to that of the surface tension: the \emph{wall free energy} is defined by
\[
\tau^{\small\mathrm{bd}}_{\beta,h}
\defby
\lim_{n\to\infty} - 
\frac{1}{(2n)^{d-1}} \log 
\frac{\sfZ^{-}_{\Delta_n;\beta,h}}{\sfZ^{\pm}_{\Delta_n;\beta,h}}.
\]
Existence of the limit is proved in~\cite{FrPf1987a}.

\subsubsection{Alternative characterizations of the wetting transition.}
We now briefly present some alternative characterizations of the wetting 
transition (a more detailed discussion can be found in the 
review~\cite{PfVe1996} and in the papers~\cite{FrPf1987a,FrPf1987b}). 
\begin{itemize}
\item \emph{Surface Gibbs states:} It can be shown~\cite{FrPf1987a,FrPf1987b} 
that complete wetting occurs 
if and only if there is a unique \emph{surface Gibbs state}, that is, if and 
only if $\lim_{n\to\infty}\mu^-_{n;\beta,h} = 
\lim_{n\to\infty}\mu^\pm_{n;\beta,h}$.
\item \emph{Thermodynamics:} It can be shown~\cite{FrPf1987a,FrPf1987b} that 
$\tau_\beta(\uvec_d) \geq \tau^{\small\mathrm{bd}}_{\beta,h}$ for 
all $\beta$ and $h$, and that partial wetting occurs if and only if 
$\tau_\beta(\uvec_d) > \tau^{\small\mathrm{bd}}_{\beta,h}$.
\item \emph{Canonical ensemble:} Here, one considers the measure 
$\mu^-_{\Delta_n;\beta,h}$ conditioned on a fixed value of the 
magnetization in $\Delta_n$: $\sum_{i\in\Delta_n} \sigma_i = m|\Delta_n|$. 
When $|m|$ is strictly smaller than the spontaneous magnetization 
$m^*(\beta)$, the system reacts by creating a droplet of $+$ phase immersed 
in a background of $-$ phase. This droplet has a deterministic macroscopic 
shape, given by the solution of a suitable generalization of the variational 
problem \textbf{(VP)} in Appendix~\ref{app:Wulff}, the so called
Winterbottom problem ~\cite{BoIoVe2000}. It can then be 
shown~\cite{PfVe1997,BoIoVe2001} that, at least on the level of the variational
problem, the resulting droplet attaches itself 
to the bottom wall if and only if \(|h|<\hw(\beta)\)
\footnote{If $h= \hw(\beta)$, then 
the Winterbottom shape 
coincides with the Wulff shape. 
It might still be attracted or repulsed by the boundary of the box. 
In particular, when $d\geq 3$, the existence of 
facets on the Wulff shape makes it possible, albeit highly unlikely,   
that the facet remains in contact with the wall. These problems are open. } 
\end{itemize}

\section{Critical prewetting}
\label{sec:crit-pre}

As we have just described, wetting phenomena occur at phase coexistence: the 
bulk of the system is occupied by some thermodynamically stable phase $A$, 
while the substrate favors another stable phase $B$, resulting, in the regime 
of complete wetting, in the creation of a mesoscopic layer of phase $B$ 
covering the substrate.

When the system is brought away from phase coexistence, phase $B$ becomes 
thermodynamically unstable and the wetting transition disappears: 
irrespectively of the preference of the substrate towards ``phase'' $B$, the 
latter never forms a mesoscopic layer. Nevertheless, there is still a trace 
of the wetting transition, at least as the system is brought close to phase 
coexistence. 
Namely, when this happens, a microscopic layer of the unstable phase $B$ 
covers the substrate, and the width of this layer diverges as the system 
approaches phase coexistence if, and only if, the system at phase coexistence 
is in the complete wetting regime.
This phenomenon is known as \emph{prewetting} and takes different forms in 
two- and three-dimensional models.

\subsection{Critical prewetting in the Ising model}

In order to discuss prewetting in the Ising model, we need to 
modify the Hamiltonian in order to bring the system away from phase 
coexistence. The most natural way to do that is by adding a positive bulk 
magnetic field, the latter ensuring that only the $+$ phase is 
thermodynamically stable. Namely, we consider the following modification of 
the Hamiltonian in~\eqref{eq:Ham:Ising2dNoBulkField}: for $\lambda\geq 0$, let
\begin{equation}
\label{eq:Ham:Ising2dWithBulkField}
\Ham_{n;h,\lambda}(\sigma)
\defby
\Ham_{n;h}(\sigma)
- \lambda \sum_{i\in\Delta_n} \sigma_i.
\end{equation}
We denote by $\mu^\pm_{\Delta_n;\beta,h,\lambda}$ the corresponding Gibbs 
measure on $\Omega^\pm_{\Delta_n}$.
Clearly, the latter reduces to $\mu^\pm_{\Delta_n;\beta,h}$ when 
$\lambda=0$. In this section, we assume that $h\geq\hw(\beta)$, that is, the 
system is in the complete wetting regime when $\lambda=0$.

\subsubsection{The two-dimensional case}

\begin{figure}[t]
\begin{center}
\includegraphics[width=\textwidth]{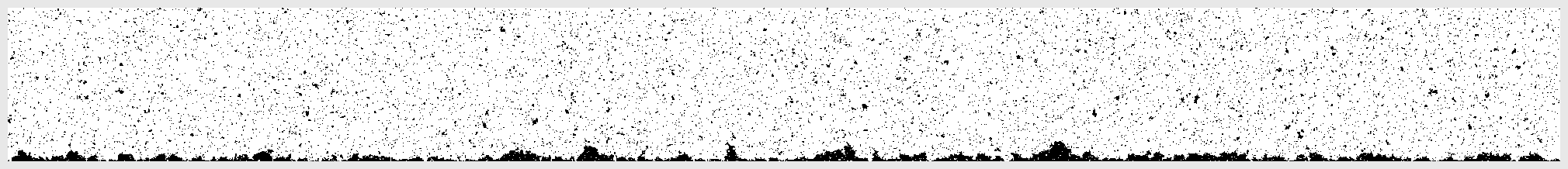}
\centerline{\small $\lambda=0.01$}

\smallskip
\includegraphics[width=\textwidth]{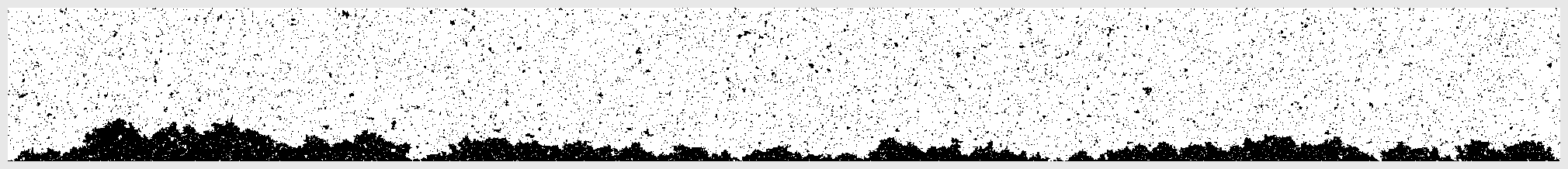}
\centerline{\small $\lambda=0.001$}

\smallskip
\includegraphics[width=\textwidth]{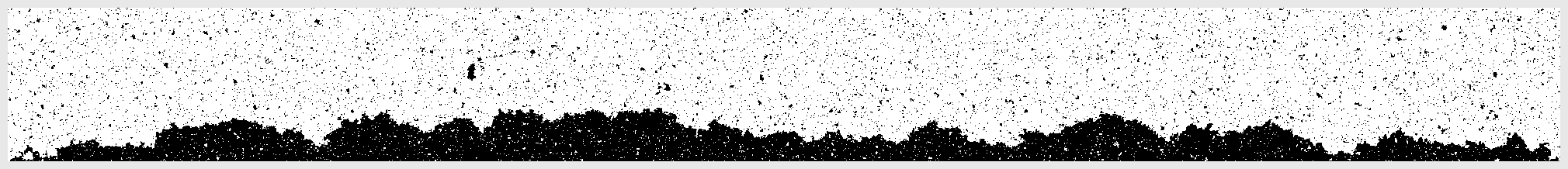}
\centerline{\small $\lambda=0.0001$}
\end{center}
\caption{Two-dimensional Ising model at $\beta=0.5$; $+$ spins are white, $-$ 
spins black. Layer of unstable $-$ phase covering the substrate in the 
complete wetting regime ($h=1$), for decreasing values of the bulk field 
$\lambda$.}
\label{fig:prewettingIsing2d}
\end{figure}%

In order to analyze prewetting in this model, we have to provide a suitable 
definition for the (average) thickness of the layer of unstable $-$ phase. 
One possible choice is to define the latter as
\[
w(\lambda;\beta,h) \defby \lim_{n\to\infty} (2n+1)^{-1}\, 
|\mathrm{below}(\gamma)|,
\]
where $|\mathrm{below}(\gamma)|$ denotes the area of the  region of $\bbR^2$ 
delimited by the unique infinite Peierls contour $\gamma$ and the line 
$\calL$. We can now state the following
\begin{conjecture}
For any $\beta>\betac(2)$ and any $h\geq\hw(\beta)$, there exist constants 
$c_1>0$, $c_2<\infty$ and $\lambda_0>0$ such that, for all 
$\lambda<\lambda_0$,
\[
c_1 \lambda^{-1/3} \leq w(\lambda;\beta,h) \leq c_2 \lambda^{-1/3}.
\]
\end{conjecture}%
We refer the reader to~\cite{Ve2004} for some partial results toward 
proving this conjecture.

The continuous divergence of the thickness of the unstable film as 
\(\lambda\downarrow 0\) is characteristic of \emph{critical} prewetting.
In fact, one has a precise conjecture for the process describing the properly 
scaled interface in the limit \(\lambda\downarrow 0\). Namely, as explained 
above, the typical width of the layer if of order \(\lambda^{-1/3}\). It is 
also possible to argue that the natural horizontal lengthscale is of order 
\(\lambda^{-2/3}\). So, one might hope that after rescaling the interface by 
\(\lambda^{1/3}\) vertically and by \(\lambda^{2/3}\) horizontally, the 
resulting process might be described by some universal scaling limit. That 
this occurs was actually established for a class of one-dimensional effective 
interface models in~\cite{IoShVe2015}. Let us describe the conjectured result
for the Ising model.

Given a particular realization of the Ising interface \(\gamma\), let us 
denote by \(\sigma_\gamma\) the configuration in \(\Omega_{\Delta_n}^\pm\) 
which has \(\gamma\) as its unique contour. We can then define the upper and 
lower ``envelopes'' \(\gamma^\pm:\bbZ\to\bbZ\) of \(\gamma\) by
\begin{gather*}
\gamma^+(i) \defby \max\setof{j\in\bbZ}{\sigma_{(i,j)}(\omega_\gamma)=-1} 
+1 ,\\
\gamma^-(i) \defby \min\setof{j\in\bbZ}{\sigma_{(i,j)}(\omega_\gamma)=+1}
-1 .
\end{gather*}
Note that \(\gamma^+(i)>\gamma^-(i)\) for all \(i\in\bbZ\).
Using the random-walk representation of Ising interfaces, described in more 
detail in Section~\ref{sec:EffectiveRandomWalk} \footnote{To be precise the effective random walk representation 
of Section~\ref{sec:EffectiveRandomWalk} with independent steps only applies at very low temperatures (large $\beta$). For moderate values 
of $\beta >\beta_c$, one needs the more refined construction of \cite{CIV03}.}, one can prove that, with 
probability close to \(1\), \(\gamma^-\) and \(\gamma^+\) remain very close 
to each other: there exists \(K=K(\beta)<\infty\) such that, with probability 
tending to \(1\) as \(n\to\infty\),
\begin{equation}\label{eq:Ising:Interface Width}
\max_{|i|\leq R} |\gamma^+(i) - \gamma^-(i)| \leq K\log R.
\end{equation}
for any \(R>0\).
Let us now introduce the rescaled profiles \(\hat\gamma^\pm:[-1,1] 
\to \bbR\). Given \(y=(y_1,y_2)\in\bbR^2\), let us write \(\lfloor y 
\rfloor \df (\lfloor y_1 \rfloor,\lfloor y_2 \rfloor)\). We then set, for any 
\(x\in\bbR\),
\[
\hat\gamma^+_\lambda(x) = \lambda^{1/3} \gamma^+(\lfloor \lambda^{-2/3}x 
\rfloor),
\]
and similarly for \(\gamma^-\). Thanks to~\eqref{eq:Ising:Interface Width}, 
for any \(c\geq 0\), 
\[
\lim_{\lambda\downarrow 0} \lim_{n\to\infty} 
\mu^\pm_{\Delta_n;\beta,h,\lambda} \bigl(\sup_{|x|\leq \lambda^{-c}} 
|\hat\gamma^+_\lambda(x) - \hat\gamma^-_\lambda(x)| < \epsilon \bigr) = 
1,\quad \forall\epsilon>0,
\]
so that in order to analyze the scaling limit of \(\gamma\), it is sufficient 
to understand the scaling limit of, say, \(\gamma^+\). 
\begin{conjecture}
\label{conj:IsingFS}
In the limit \(n\to\infty\), followed by \(\lambda\downarrow 0\), the 
distribution of \(\hat\gamma^+_\lambda\) under 
\(\mu^\pm_{\Delta_n;\beta,h,\lambda}\) converges weakly to that of a 
Ferrari--Spohn diffusion with parameters \(\sigma^2 = \chi_\beta \) 
and \(q(r)=2m^*r\) 
(see Appendix~\ref{app:FS}), where $\chi_\beta = \lb\tau_\beta (0) + 
\tau_\beta^{\prime\prime} (0)\rb^{-1}$ is the curvature 
of the unnormalized equilibrium crystal shape $\calW_\beta$ in lattice 
directions $\pm \sfe_i$ 
 and 
\(m^*=m^*(\beta)\df \lim_{n\to\infty} 
\mu^+_{\Lambda_n;\beta}(\sigma_0)\) is the spontaneous magnetization.
\end{conjecture}
Recall that Wulff shape are defined and discussed in  Appendix~\ref{app:Wulff}.
In  Appendix~\ref{app:SAW} we try to justify the above conjecture by sketching
an argument for proving a similar statement in the context of self-avoiding 
walks under area tilts.

\subsubsection{The three-dimensional case.}

As mentioned above, infinitely many layering transitions are conjectured to 
occur in the partial wetting regime in this case. Not surprisingly, this has 
also consequences for the prewetting behavior: 
namely, as the bulk field $\lambda$ is decreased towards zero, the continuous 
divergence observed in the two-dimensional model is replaced by an infinite 
sequence of first-order phase transitions, also known as layering 
transitions. The latter are however much better understood rigorously than 
the corresponding transitions in the absence of a bulk external field.

The existence of the first such layering transition was established, at low 
temperatures, in~\cite{FrPf1987b}.

Moreover, in the case $h=1$, it was argued
\footnote{The claim below 
 is formulated for the Ising model proper and it is 
is much stronger than the corresponding statements 
in \cite{DiMa94,CeMa96}, which hold for the SOS simplification - 
see Subsection~\ref{subsub:SOS}.  
We have not checked the proof in \cite{Ba07}
}
in~\cite{Ba88,Ba07} 
that, there is 
a decreasing sequence $(\lambda(n))_{n\geq 0}$ such that $\lim_{n\to\infty} 
\lambda(n) = 0$ and
\begin{itemize}
\item for each $n\geq 0$, when $\lambda \in (\lambda(n),\lambda(n+1))$, the 
film of unstable $-$ phase has a thickness of $n$ microscopic layers (with a 
microscopically sharply-defined boundary);
\item at the transition points $\lambda=\lambda(n)$, there is a coexistence 
of two Gibbs states with layers of thickness $n-1$ and $n$.
\end{itemize}
In addition, some information on the location of the first few lines can be 
found in these papers.
In principle, delocalization and scaling of these top level lines as the bulk 
field $\lambda$ goes to zero, should bear resemblance to critical
prewetting of Ising interfaces in two dimensions. 

\subsubsection{Relation to metastability}

In this section, we describe a different setup in which a closely related
phenomenon occurs. We consider only the two-dimensional Ising model 
in the box $\Lambda_n = \{-n, \dots , n\}^2
$ with \(-\) boundary condition and the 
Hamiltonian~\eqref{eq:Ham:Ising2dWithBulkField}, but something similar occurs 
also in higher dimensions.
We will be interested in the behavior of this model when $h=1$ and 
$\lambda>0$.
That is, we consider the Gibbs measure
\(
\mu^-_{\Lambda_n;\beta,1,\lambda}
\).

As before, because of the presence of a positive bulk field $\lambda$, the 
$+$ phase is the unique equilibrium phase. Now, however, the boundary 
condition favors $-$ spins, which leads to delicate metastability issues.
Namely, two types of behavior could be expected: either the boundary 
condition dominates and the box is filled with the \(-\) phase, or the bulk 
field dominates and the box is occupied by the \(+\) phase (except, possibly, 
close to the boundary). Replacing the \(-\) phase inside the box by the \(+\) 
phase yields an energetic gain of order \(\lambda |\Lambda_n| \) due to 
preference of the bulk field for the \(+\) phase, but there is also an 
energetic cost of order \(|\partial\Lambda_n|\) associated to the boundary of 
the ``droplet'' of \(+\) phase thus created. One would thus expect that the 
transition between these two regimes should occur for a value of \(\lambda\) 
of order \(1/n\). This is indeed the case.

\begin{figure}[t]
\begin{center}
\includegraphics[width=.5\textwidth]{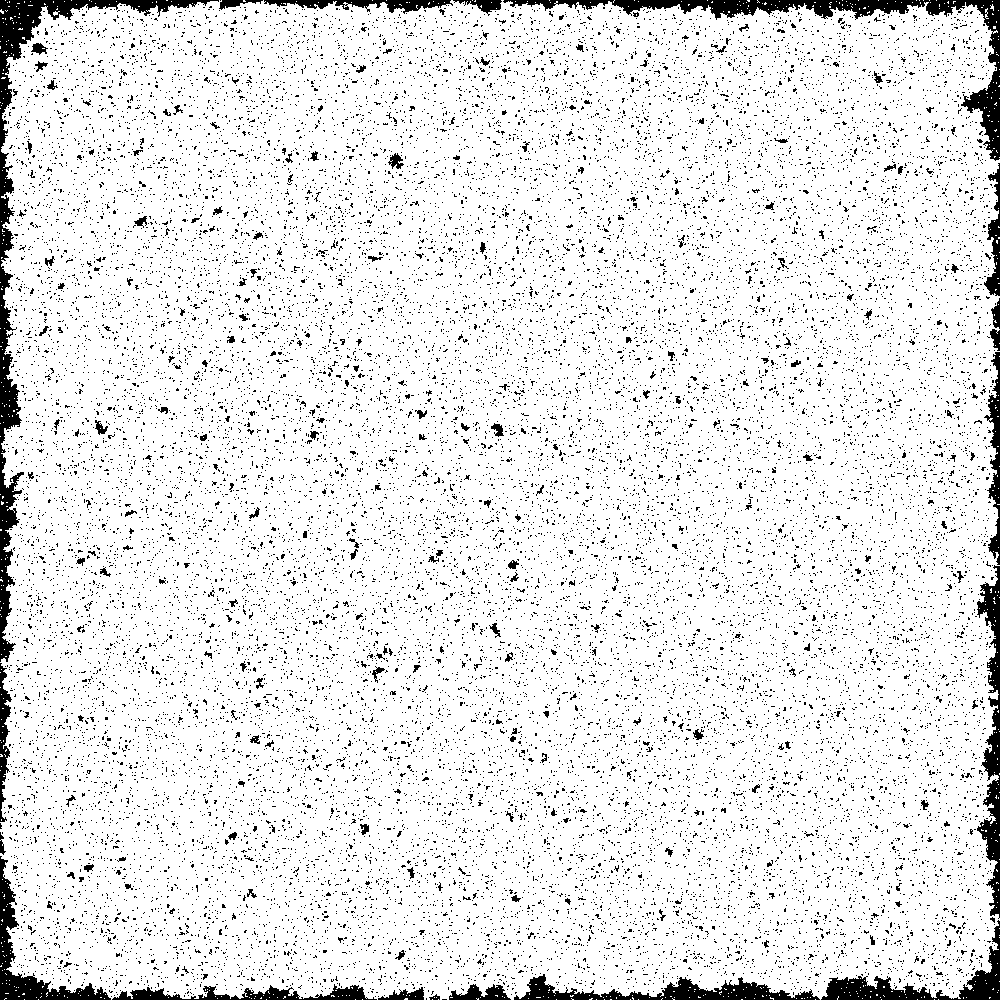}
\end{center}
\caption{The two-dimensional Ising model at $\beta=0.48$ in $\Lambda_{1000}$ 
with $-$ boundary condition and a bulk magnetic field $\lambda=0.007$.}
\label{fig:MetastabilityIsing}
\end{figure}

More precisely, the following is proved in~\cite{SS96}:  Let $\sfW_\beta$ be 
the unit area Wulff shape for the surface tension $\tau_\beta$, 
precisely as specified in Appendix~\ref{app:Wulff}. Let $\nu_\beta$ be 
the critical slope for the dual constrained 
variational problem {\bf (DCVP)}, see formula 
\eqref{eq:crit-beta}. 
Let 
\be{eq:B-not}
B_0(\beta)\defby 
\frac{1}{2m^* (\beta )} \nu_\beta = 
\frac{4\tau_\beta(\sfe ) + \tau_\beta(\sfW_\beta)}{8 m^*_\beta} .
\ee
Then, for all \(\beta>\betac(2)\), there exists \(K(\beta)<\infty\) such 
that the following 
holds:
\begin{itemize}
\item Let \(B<B_0(\beta)\) and set \(\lambda = \lambda(n) = B/n\). Then, as 
\(n\to\infty\), with probability going to \(1\), all Peierls contours are of 
diameter at most 
\(K(\beta)\log n\)~\footnote{Strictly speaking in \cite{SS96} it is only 
proven that the small contours are at most of length 
\(\epsilon/\lambda=\epsilon n/B\), but 
the claimed extension should be routine using the techniques introduced in
\cite{ISc98}.}.
\item Let \(B>B_0(\beta)\) and set \(\lambda = \lambda(n) = B/n\).
Define $\nu = 2m^* (\beta ) B > \nu_\beta$, and let $\sfa = \sfa (\beta , \nu )$
 be the corresponding solution of \textbf{(DCVP)}, 
 and $\sfP^\sfa_\beta$ be the corresponding Wulff plaquette. \newline 
 Then, for any \(\epsilon>0\), the following occurs with probability going to 
\(1\) as \(n\to\infty\) (see Figure~\ref{fig:MetastabilityIsing}): 
 There is a 
unique Peierls contour of diameter more than \(K(\beta)\log n\) and this contour is 
contained in the set \((1+\epsilon)n \sfP^a_\beta \) and surrounds the set 
\((1-\epsilon)n\sfP^a_\beta \). 
\end{itemize}

Therefore, in the regime \(B>B_0(\beta)\), there are macroscopic interfaces 
along (part of) the four walls. As before, these interfaces coincide with the 
walls only at the macroscopic scale. At the microscopic scale, there is a 
layer of \(-\) phase between the walls and the droplet of \(+\) phase. It is 
conjectured that the width of this layer is of order
\[
\lambda(n)^{1/3} \sim n^{-1/3},
\] 
as \(n\to\infty\), in perfect analogy with the discussion of prewetting above.
Moreover, the $(n^{2/3} , n^{1/3})$-rescaling of Ising interfaces along flat
boundaries of the box is expected to lead to 
 Ferrari--Spohn asymptotics as described in Conjecture~\ref{conj:IsingFS}.

\subsection{Interfacial wetting and prewetting in the 2d Blume--Capel model}

The phenomena described above in the two-dimensional Ising model occur of course
in a broad class of systems. As one example, we discuss here the two-dimensional
Blume--Capel model. The latter has spins taking values in the set \{-1,0,1\} and
a formal Hamiltonian of the form
\[
-J \sum_{i\sim j} \sigma_i\sigma_j - \lambda \sum_{i} \sigma_i^2 .
\]
It is well known~\cite{BrSl} that, at sufficiently low temperatures, there exists
\(\lambda_c=\lambda_c(\beta)>0\) such that, for all \(\lambda<\lambda_c\), there
exists a unique extremal translation-invariant Gibbs state, typical
configurations of which are composed of an infinite ``sea'' of \(0\) spins with
only finite ``islands''. In contrast, for all \(\lambda>\lambda_c\), there are
exactly two extremal translation-invariant Gibbs states with ``seas'' of \(+\),
respectively \(-\) spins. Finally, at the triple point \(\lambda_c\), all three
measures coexist (and are the only extremal translation-invariant states).

Let us consider the model in the box \(\Lambda_n\) with a boundary condition made
up of \(+\) spins in the upper half plane, and \(-\) spins in the lower half
plane. Assume first that \(\lambda=\lambda_c\). Then, as in the Ising model, this
boundary condition favors \(+\) spins in the top half of the box and \(-\) spins
in the bottom half, with an interface separating the two regions. However, now a
new phenomenon occurs. Let us denote by \(\tau^{+-}\) the surface tension
associated to the interface between the \(+\) and \(-\) phases, and define
similarly \(\tau^{+0}\) and \(\tau^{-0}\). It is expected that
\(\tau^{+-}=\tau^{+0}+\tau^{-0}\) and that this relation should imply that it is
preferable for the system to introduce a mesoscopic layer of \(0\) phase between
the regions occupied by the \(+\) and \(-\) spins. The width of this layer
should be of order \(\sqrt{n}\). No rigorous results seem to have been obtained
on this problem, but a discussion can be found in~\cite{BrLe-BlumeCapel}.

Let us now assume that \(\lambda>\lambda_c\). In this case, the 
layer of \(0\) phase becomes unstable and, similarly to what we saw 
for prewetting, it becomes microscopic. One can then once more 
analyze how the width of this layer changes as \(\lambda\downarrow\lambda_c\). 
The natural conjecture, supported by numerical simulations (see~\cite{SeYe} for 
the earliest we could find), is that it should again behave as \(\lambda^{1/3}\), 
which suggests a Ferrari--Spohn structure of appropriate scaling limits. 
Again, none of the above has 
yet been established, although some related results have been 
obtained far away from the triple point, \(\lambda \gg \lambda_c\)~\cite{HrKo}.

\section{Effective interface models} 
\label{sec:Effective}

The discrete two-dimensional effective interface models discussed in this 
section are defined over lattice boxes  $\Lambda_N\subset \bbZ^2$,
\be{eq:Lambda-N}
\Lambda_N = \{-N,-N+1, \dots, N-1, N\}^2,
\ee
$N$ being the linear size of the system. The interface is described as 
the set of \emph{integer} random heights $\varphi  = 
\{\varphi_x\}_{x\in\bbZ^2}$. Thus, $\varphi_0$ describes the random height of 
the interface at the center of the box $\Lambda_N$. 

Unless mentioned otherwise we shall assume that the interfaces in question 
are pinned at zero height  outside $\Lambda_N$; $\varphi_x = 0$ if 
$x\not\in\Lambda_N$. 

The passage from the Ising model with boundary conditions~\eqref{eq:pm-bc} 
to the effective interface models considered in this Section can be viewed 
either as a strong interaction limit, in which the coupling constants 
$J_{i,j}$ in~\eqref{eq:Ham:Ising2dNoBulkField} are set to \(+\infty\) for 
vertical  bonds $(i, j)$, or as a simplification of ``true'' Ising interfaces 
obtained by prohibiting overhangs. The former point of view gives rise to the 
class of Solid-on-Solid (SOS) models specified in Subsection~\ref{subsub:SOS} 
below. The latter point of view still leaves room for incorporating 
low-temperature phases above and below the interface. An effective model
of this sort is discussed in Subsection~\ref{subsub:SOS-bulk}.

\subsection{Low-temperature interfaces and ensembles of level lines} 
At low temperatures, the large-scale properties of effective interfaces 
should be recovered from the statistics of large level lines which locally 
look like one-dimensional effective random walks. 

\subsubsection{SOS-type models}  
\label{subsub:SOS}
The Hamiltonian $\calH_N$ and the probability distribution $\bbP_N^\beta$ are
given by
\be{eq:HamPM-SOS}
 -\calH_N (\varphi )
 =
 \sum_{xy}J_{xy} U (\varphi_x - \varphi_y ) 
 +
 \sum_x f (\varphi_x )
 \quad \text{ and }\quad 
 \bbP_N^\beta (\uphi )
 =
 \frac{1}{Z_N } {\rm e}^{-\beta \calH_N (\uphi )}, 
\ee
where $Z_N = Z_N (\beta, V, J, f )$ is the partition function. Familiar 
examples in the case of nearest-neighbor interactions, $J_{xy}=\1_{x\sim y}$, 
and identically zero self-potential, $f\equiv 0$, include: 

\begin{enumerate}[label=(\alph*)]
\item\label{exSOS}
SOS model if $U(t ) = \abs{t}$.
\item
Discrete Gaussian (DG) surface if $U(t ) = |t|^2$.
\item
Restricted SOS model in the case of $U (t) = \abs{t} +\infty\cdot\1_{\abs{t} 
>1}$. 
\end{enumerate}

Below, we shall restrict attention to~\ref{exSOS}, but other potentials 
should be kept in mind. 

The self-interaction $f$ is used to model external influences, such as bulk 
fields, interaction with 
a substrate or a hard-wall constraint. 
Consider, for instance, 
\be{eq:Self-f} 
f(\varphi)
=
-\infty\cdot\IF{\varphi <0} - \lambda \varphi +  h\,\IF{\varphi = 0} . 
\ee
The first term forces the surface to stay above the wall: $\varphi_x \geq 0\; 
\forall x$. 
If we use~\eqref{eq:HamPM-SOS} to model a random interface between coexisting 
phases, then the parameter $\lambda$ measures the strength of the bulk 
magnetic field, while $h$ measures the interaction between the interface and 
a substrate placed at zero height. 

\smallskip
Random interfaces modeled with Hamiltonians as in~\eqref{eq:HamPM-SOS} 
display a rich phenomenology. 
\smallskip 

\paragraph{\bf Roughening.}
Let $f\equiv 0$. It is well known (and not difficult to prove) that 
\(\bbE_N^\beta {\abs{\varphi_0}}\) is bounded uniformly in \(N\) when \(\beta\) is 
large, for both the SOS and DG models. Groundbreaking results by Fr\"{o}hlich 
and Spencer~\cite{FS81} imply that there is a roughening transition for both 
models. Namely, 
\(\lim_{N\to\infty}\bbE_N^\beta {\abs{\varphi_x-\varphi_0}^2} \geq c\log |x|\). 
The result is proved by a mapping to a Coulomb gas; finding a direct 
probabilistic argument and deriving a (presumably Gaussian) scaling limit are 
well-known challenges. 
\smallskip 

\paragraph{\bf Entropic repulsion.}
Let $f(\varphi) = -\infty\cdot\IF{\varphi<0}$. 
Early results~\cite{BMF86, LM87} imply that the interface, as the size $N$ of 
the system grows, is always repelled to infinity, even in the case of large 
$\beta$ for which the unconstrained interface is localized. 

In the context of SOS models, the lower temperature (that is, $\beta \gg 1$)
situation was analyzed in much more detail in the recent 
works~\cite{CLMST14, CLMST13, CMT14}. 
Their results include sharp height concentration (on a 
diverging level $\sim\log N$), macroscopic scaling limit for level sets, and 
$N^{1/3}$ bounds on fluctuations of level lines long the flat segments of 
the boundary $\partial\Lambda_N$.  
Furthermore, height concentration results (but not macroscopic scaling limits 
and approximate order of fluctuations of level lines) were established for 
low-temperature models with more general interactions of the form $U(t) = 
\abs{t}^p\; p\in (1,\infty]$ in~\cite{LMS14}.
\smallskip

\paragraph{\bf Prewetting and Layering.} 
Let $f(\varphi) = -\lambda \varphi - \infty\cdot\IF{\varphi<0}$. 
The bulk field $\lambda$ penalizes large values of the volume below the 
surface; in this sense, it models the bulk magnetic field 
in~\eqref{eq:Ham:Ising2dWithBulkField}.
When $\lambda>0$, there is a competition between the entropic repulsion and 
the volume penalization imposed by $\lambda$. In the \emph{low-temperature} 
regime, this results is an infinite sequence of first-order phase 
transitions, which can be characterized in terms of the concentration of the 
typical interface height at increasing values: $1,2,3,\dots$; 
see~\cite{DiMa94, CeMa96, LMazel96} for precise statements of the results. 
Roughly speaking, it was proved that for any sufficiently large $\beta$ there 
exists a number $n_{{\rm max}} (\beta )$ and a sequence of bulk fields 
$\lambda_0 (\beta ) >\lambda_1 (\beta ) > \dots >\lambda_{n_{\rm max}}(\beta )$, 
such that  the surface is concentrated on height $n$ whenever 
\be{eq:SOS-heights-n} 
\lambda \in \lb \lambda_{n} (\beta ), \lambda_{n-1} (\beta )\rb . 
\ee
It is not clear whether the restriction $n\leq n_{{\rm max}}$ is of a technical 
nature or not. Presumably it is not, and 
prewetting thus occurs as $\lambda\downarrow 0$ through an infinite sequence of
jumps of the typical interface height. 

\smallskip 

\paragraph{\bf Wetting.}
Let $f(\varphi) = h\, \IF{\varphi=0} - \infty\cdot\IF{\varphi <0}$. The 
boundary field $h$ plays the same role as in~\eqref{eq:J-h}. There is now a 
competition between entropic repulsion and attraction by the substrate 
(measured through \(h\)).  Define the {\em dry} set $\calA_N = 
\setof{x\in\Lambda_N}{\varphi_x=0}$. For SOS interfaces, it was proved 
in~\cite{Chalker} that typically $|\calA_N| \sim N^2$ when $h $ is large, 
while $|\calA_N| \ll N^2$, 
when $h>0$ is small.
More precisely, it is proven in~\cite{Chalker} that the dry set has a 
uniformly (in the linear size of the system $N$) positive density if
\be{eq:Chalker1} 
\beta h
>
-\log\Bigl(\frac{1 -{\rm e}^{-\beta}}{16 (1 + {\rm e}^{\beta})}\Bigr) ,
\ee
whereas the density of the dry set vanishes
($|\calA_N|$ being at most of order $N$), 
as $N\to\infty$, as soon as 
\be{eq:Chalker2} 
\beta h < - \log (1 -{\rm e}^{-4\beta}).
\ee
This is a perturbative result, which {does not directly address} the
nature of transition
for $h$-s between the two values in~\eqref{eq:Chalker1}
and~\eqref{eq:Chalker2}.

The \emph{lower temperature} situation was analyzed in much more detail 
in~\cite{AlDuMi11}, where the existence of a sequence of layering transitions 
has been established and explained. Roughly speaking, the following is shown 
(see Theorem~1.1 in~\cite{AlDuMi11} for the precise statement). Fix
$\epsilon>0$. Then, for any given $n\in \bbN$ and for all $\beta$ sufficiently
large, the surface concentrates on height $n$ whenever
\be{eq:Al-wet-cond} 
-\ln\lb 1 - {\rm e}^{-4\beta}\rb + (2+\epsilon ){\rm e}^{-2\beta (n+3)} \leq 
h\beta \leq -\ln\lb 1 - {\rm e}^{-4\beta}\rb + (2-\epsilon ){\rm e}^{-2\beta (n+2)} .
\ee
Such a result had previously been stated, but
not proved, in~\cite{BMF86}. 

Finally, we note a  recent preprint \cite{Lacoin2017a} 
which contains (Theorem~2.1 there) an identification of the right hand side of \eqref{eq:Chalker2} as the 
low temperature wetting transition threshold and, furthermore, contains a claim on an almost piece-wise
affine structure which bears a strong indication of an infinite sequence of layering transitions. 
A proof and a comprehensive analysis 
of such infinite sequence of layering transitions appeared (Theorems 2.1 and 2.6 there) in a very recent
sequel \cite{Lacoin2017b} to \cite{Lacoin2017a}. We have not checked the arguments of 
\cite{Lacoin2017a, Lacoin2017b} in detail. If correct (which is presumably the case) they fetch a 
complete, modulo a certain analyticity 
issue - see Remark~2.2. in \cite{Lacoin2017b}, original description of layering transition driven by boundary
fields in the SOS context. 

\subsubsection{SOS-type model coupled to a bulk field}
\label{subsub:SOS-bulk}

The interfaces have the same distribution as in~\eqref{eq:HamPM-SOS}, with 
$f\equiv 0$. Instead of the latter, consider a coupling with high- and 
low-density bulk Bernoulli fields below and, respectively, above 
the interface. These fields  are designed to mimic co-existing low- and high- 
density phases, say vapor and solid, in lattice gases. Namely, consider a 
three-dimensional vessel 
\[
B_N = \Lambda_N\times \{ -N +\tfrac{1}{2}, \ldots , N-\tfrac{1}{2} \}.
\]
A realization of $\varphi$ (constrained to stay inside $B_N$)
splits $B_N$ into two halves; $B_N = S_N(\uphi)\cup V_N(\uphi)$. 
Given $0 < p_{\rm v} < p_{{\rm s}} < 1$, we place particles with probability 
$p_{{\rm s}}$ into sites of $S_N$ and with probability $p_{\rm v}$ into sites 
of $V_N$. We use $\bbQ_N$ to denote the distribution of the resulting triple 
$(\varphi, \eta^{{\rm s}} ,\eta^{{\rm v}})$. That is,
\be{eq:B-SOS}
\bbQ_{N}^\beta (\varphi, \eta^{{\rm s}}, \eta^{{\rm v}})
\propto {\rm e}^{-\beta \calH_N(\varphi)} 
\bbB_{S_N(\uphi)}^{p_{{\rm s}}}(\eta^{{\rm s}})  
\bbB_{V_N(\uphi)}^{p_{{\rm v}}}(\eta^{{\rm v}}),  
\ee
where $\bbB^p_A$ is the $p$-Bernoulli product measure on $\{0,1\}^A$.  

\smallskip
\paragraph{\bf Faceting.} 
The measure~\eqref{eq:B-SOS} does not incorporate a bulk chemical potential, 
a hard wall or the interaction with an active substrate. 
However, the layering phenomenon can be modeled in the following way:
Note that the average number of particles under $\bbQ_N^\beta$ is 
$(p_{{\rm s}} + p_{{\rm v}})N \abs{\Lambda_N} \df r_N N^3$.
At \emph{low temperatures}, the interface is flat. 
Facets of macroscopic size appear under the canonical constraint
\be{eq:BSOS-Q}
\bbQ_{N,\sfa}^\beta (\cdot) =  
\bbQ_{N}^\beta \bigl(\,\cdot \bigm| \sum_{\sfz\in S_N} 
\eta^{\rm s}_\sfz + \sum_{\sfz\in V_N}\eta^{{\rm v}} = 
\lfloor r_N N^3 + \sfa N^2\rfloor \bigr) .
\ee
The model in~\eqref{eq:B-SOS}, for the SOS interaction $U(t) = \abs{t}$,
was introduced in~\cite{IS08}, where it was established that first-order 
phase transitions under the canonical measure~\eqref{eq:BSOS-Q} occur (as the 
parameter $a$ grows), corresponding to the creation of the first two facets. 
Actually the forthcoming work~\cite{IS16} implies that the model undergoes 
an infinite sequence of such first-order layering transitions. 
The discontinuity of facet sizes at critical values of $a$ is due to the 
coupling with the bulk Bernoulli field, in a way similar to the spontaneous 
creation of a droplet in the 2D Ising model~\cite{BCK03}.
Presumably, coupling of an SOS interface with bulk Bernoulli fields reflects 
in a more accurate way the faceting of microscopic Wulff crystals in the
low-temperature 3D Ising model, which was the motivation for an earlier study 
of the phenomenon in the context of pure SOS models in~\cite{BSS05}. 

\subsubsection{Level lines of low temperature interfaces} 
\label{susub:level-lines}

Random interfaces under either~\eqref{eq:HamPM-SOS} or~\eqref{eq:B-SOS} can 
be described in terms of their level lines. Since  the surface is pinned  at 
zero height outside $\Lambda_N$, level lines are closed microscopic contours. 
It is convenient to employ north-east splitting rules to avoid ambiguities at 
vertices which are incident to four edges, see for instance 
\cite[Subsection~2.1]{IofShlosTon15} or~\cite[Definition~2.1]{CLMST13}. In 
the sequel we shall tacitly assume that such rules are employed and, 
consequently, we shall talk about \emph{compatible} portions of level lines 
or contours without further comments. 

Given a contour $\gamma$, we use $\abs{\gamma}$ to denote its length, 
$\mathrm{diam}(\gamma)$ to denote its diameter in $\|\cdot\|_\infty$-norm, 
$\mathring{\gamma}$ to denote its interior and $\sfA = \sfA(\gamma)$ to 
denote the area of $\mathring{\gamma}$.
Fix $\epsilon > 0$ sufficiently small. 
If $N$ is the linear size of the system, a microscopic contour $\gamma$ is 
said to be large if $\mathrm{diam}(\gamma) \geq \epsilon N$ and
small if $\mathrm{diam}(\gamma) \leq \epsilon^{-1}\log N$. 
Otherwise, it is said to be intermediate. 

In the case of low-temperature SOS ($U(t) = \abs{t}$) Hamiltonians, contour 
weights factorize and cluster expansions were derived and explored, 
see for instance~\cite{DKS, 
DiMa94, CeMa96, Miracle99, IS08, AlDuMi11, CLMST14, CLMST13, IofShlosTon15}. 
It is natural to make expansions relative to the unconstrained SOS model 
(with $f\equiv 0$). In all situations of interest, one should \emph{in principle} 
be able to rule  out intermediate contours 
\footnote{At least on a heuristic level. Rigorous results were derived and formulated
differently in different works; for instance, in~\cite{AlDuMi11} the authors 
do not fix the interface at zero height outside of a finite box, but rather 
directly explore stability properties of infinite-volume interfaces which live on 
heights $n=0, 1, 2, \dots $}
and check that there is at most one ordered stack 
\be{eq:contour-stack}
\ugamma = \{\gamma_1,\ldots,\gamma_n\} \text{ satisfying }  
\mathring{\gamma}_1\subseteq \mathring{\gamma}_2\subseteq\cdots\subseteq 
\mathring{\gamma}_n \subseteq \Lambda_N
\ee
of compatible large microscopic contours. 
A somewhat simplified version of the SOS polymer weights 
$q_{\beta,N}(\ugamma)$ that are obtained via cluster expansions with 
respect to small contours reads as:
\be{eq:Form-Weights}
\log q_{\beta,N}(\ugamma)
=
-\beta \sum \abs{\gamma_i }
+
\sum_i\sum_{\calC} \Phi_{\beta,N}(\calC;\gamma_i)
+
\Psi_{\beta,N}(\ugamma) 
+
F_{\beta,N}(\sfA_1,\sfA_2,\ldots) .
\ee
Above $\Phi_{\beta,N}(\calC;\gamma)$ are the cluster weights for the model 
defined on a finite box $\Lambda_N$; in particular, $\Phi_{\beta,N}$ 
incorporates the effects due to the finite geometry of the latter, 
including short-range (exponential tails) interactions with the boundary 
$\partial\Lambda_N$. 
$\Psi_{N,\beta}$ is a multi-body short-range (again in the sense of 
exponential decay with $\mathrm{diam}(\calC)$) interaction between the 
different contours in a stack, which reflects the possibility of sharing 
clusters. $F_{\beta,N}$ takes different forms in different models. Let us 
provide some examples.

\smallskip
\paragraph{\bf SOS model coupled to a bulk Bernoulli field}
Under~\eqref{eq:BSOS-Q}, the function $F_{\beta,N}$ takes,
up to $\log N$-corrections, the following form:
\be{eq:F-BSOS}
F_{\beta,N}(\sfA_1,\sfA_2,\ldots) 
=
-\frac{(\sfa N^2 -\sum_i \sfA_i )^2}{2 D \abs{B_N}}, 
\ee
where $D = \tfrac12\{p_{\rm s} (1- p_{\rm s}) + p_{\rm v} (1- p_{\rm v})\}$. 
It should be clear where~\eqref{eq:F-BSOS} is coming from: The volume of the 
surface $\varphi$ with the prescribed collection $\ugamma$ of large 
microscopic contours should concentrate around $\sfA_1+\dots + \sfA_n$.
The excess number of particles is $\sfa N^2 - \sum_1^n \sfA_i$, as imposed by 
the canonical constraint in~\eqref{eq:BSOS-Q}, and it should be compensated 
by fluctuations of the density of the Bernoulli fields above and below the 
surface. The latter have an average variance $D$ per site. 

\smallskip
For the Hamiltonians~\eqref{eq:HamPM-SOS}, \eqref{eq:Self-f}, the 
$\bbZ$-symmetry of infinite-volume states is broken and there are different
values of metastable free energies, which one computes using restricted contour
ensembles, $f^*_m= f^*_m (\beta , \lambda , h)$  for interfaces living on 
heights $m\in \bbN_0$.  Accordingly, the expression for  $F_{\beta,N}$ should
take into account these metastable free energies. Let us write it down 
up to lower-order terms for surfaces which are represented by ordered
stacks~\eqref{eq:contour-stack}. For such a surface, the total area of the
interface at height $m=0, \dots , n$
\be{eq:araeas} 
\begin{cases}
 \sfA_1, \qquad\qquad &\text{if $m=n$,}\\
 \sfA_2 - \sfA_1 &\text{if $m= n-1$,}\\
 \vdots &\,  \\
 \abs{\Lambda_N} - \sfA_n &\text{if $m = 0.$} 
\end{cases}
\ee
Hence, 
\begin{align}
\label{eq:FN-form-SOS} 
F_{\beta , N}  (\sfA_1, \sfA_2, &\dots, \sfA_n )\notag\\
&= 
f^*_n \sfA_1 + f^*_{n-1}(\sfA_2 - \sfA_1 )  + f^*_{n-2}(\sfA_3 - \sfA_2 ) +\dots + 
f_0^* (\abs{\Lambda_N } - \sfA_n )\notag\\
&=
\sum_{\ell = 1}^{n} ( f_{n-\ell +1}^* - f^*_{n-\ell}) \sfA_\ell
+ f^*_0 \abs{\Lambda_N } 
\df
\sum_{\ell = 1}^{n} r^*_{n-\ell+1} \sfA_\ell + f^*_0 \abs{\Lambda_N } . 
\end{align}

\smallskip
\paragraph{\bf Entropic repulsion}
In the case of a pure hard-wall constraint,
$f(\varphi) = -\infty\cdot \IF{\varphi < 0}$, it was established
in~\cite{CLMST13} (see~(1.6) in the latter paper) that (again up to lower-order
terms) free energies $f_m^*$ and, accordingly, area tilts $r^*_m$ are of the
form
\be{eq:SOS-lop}
f_m^* = - c_\beta {\rm e}^{-4\beta (m +1)}\quad {\rm and}\quad 
r^*_m = c_\beta {\rm e}^{-4\beta m}\lb 1 - {\rm e}^{-4\beta}\rb.
\ee
Consequently, 
\be{eq:F-SOS} 
F_{\beta,N}(\sfA_1,\sfA_2,\ldots,\sfA_n) 
= 
\bar{c}_\beta 
\sum_{\ell = 1}^{n} {\rm e}^{ -4\beta (n-\ell +1 )} A_\ell + f^*_0 \abs{\Lambda_N } .
\ee
where we have set $\bar{c}_\beta = c_\beta \lb 1  -{\rm e}^{-4\beta}\rb$. 
By Lemma~2.4 in \cite{CLMST13} the limit $\lim_{\beta\to\infty} \bar{c}_\beta =1$. 

\smallskip
\paragraph{\bf Bulk magnetic field.} In the case of pure entropic repulsion,
the area tilts $r^*_\ell$ in~\eqref{eq:SOS-lop} are always positive. 
In the presence of bulk field,  $f(\varphi) = -\lambda\,\varphi 
-\infty\cdot \IF{\varphi < 0}$, the situation is different. Roughly speaking,
it is shown in~\cite{CeMa96} (see Theorem~4.1 there for the  precise statement) 
that $r^*_n$ is positive as long as $\lambda < \lambda_{n-1}(\beta)$, 
see~\eqref{eq:SOS-heights-n}. This means that, on the level of resolution suggested
by~\eqref{eq:Form-Weights} there could be at most 
\be{bar-n-balk} 
\bar{n} = \bar{n} (\beta , \lambda ) = \max\setof{n}{\lambda \leq \lambda_{n-1} (\beta )}
\ee
large contours in the stack \eqref{eq:contour-stack}, and, accordingly, 
we should restrict attention to
\be{eq:F-SOS-B}
F_{\beta,N}(\sfA_1,\sfA_2,\ldots)  = 
\sum_1^{\bar{n}} r_{\bar{n} - \ell +1}  
\sfA_\ell  +  f^*_0 \abs{\Lambda_N} .
\ee

\smallskip 
\paragraph{\bf Boundary magnetic fields.} 
The case $f(\varphi) = -\infty\cdot \IF{\varphi < 0} + h\1_{\varphi = 0}$
is even more complicated. Indeed, \eqref{eq:Al-wet-cond} is based on 
Proposition~3.2 in~\cite{AlDuMi11}, which indicates that the area tilt $r^*_n$ is
positive as long as 
\be{eq:n-wet-cond} 
\beta h
\leq
-\ln\lb 1 - {\rm e}^{-4\beta}\rb + (2-\epsilon ){\rm e}^{-2\beta (n+2)} .
\ee
This means that, as in the case of bulk field, one should restrict attention to
a finite number $\bar{n}(\beta,h)$ of large contours in the
stack~\eqref{eq:contour-stack}.

\subsection{Variational problems and macroscopic scaling limits.} 
\label{sub:var} 

In order to discuss the macroscopic scaling limits of large microscopic level
lines for the $2+1$ dimensional surfaces we consider here, we need to 
introduce the notion of \emph{two-dimensional inverse correlation length}
\footnote{Strictly speaking, the name \emph{surface tension} might be more 
appropriate; we decided to use \emph{inverse correlation length} in order to
stress the difference with the \emph{general $d$-dimensional} Ising surface tension 
$\tau_\beta$ introduced in Subsection~\ref{sub-terminology}. Of course, if we 
are talking about interfaces in two-dimensional Ising model at sufficiently low
temperatures, then $\xi_\beta = \tau_\beta$}, as well as of Wulff shapes 
and Wulff plaquettes, which are related to the low-temperature massive 
structure of these level lines. 

\subsubsection{Inverse correlation length}  
For low-temperature polymers with infinite-volume weights
\be{eq:polym-weights}
\log q_\beta^{\sf f}(\gamma)
=
-\beta \abs{\gamma}
+
\sum_{\calC} \Phi_{\beta}(\calC;\gamma),
\ee 
the inverse correlation length $\xi_\beta$ is well defined; see, e.g., 
\cite{IS08,CLMST13} or~\cite[Subsection~2.2]{IofShlosTon15}.
Specifically, for $\sfx\in\bbZ^2$, set
$q_\beta^{\sf f}(\sfx)
= 
\sum_{\gamma: 0\to\sfx} q_\beta^{\sf f}(\gamma)$.
Then, given a direction $\frn\in\bbS^1$, 
\be{eq:tau-beta} 
\xi_\beta(\frn)
=
-\lim_{M\to\infty} \frac{1}{M}\log q_\beta^{\sf f}(\lfloor M\frn\rfloor).   
\ee
Finally, $\xi_\beta$ extends to $\bbR^2$ by homogeneity. 
More details on the existence and properties  of the limit 
in~\eqref{eq:tau-beta} can be found in Subsection~\ref{sub:Ising} below. 
Note that, for the models we consider, the cluster weights 
in~\eqref{eq:polym-weights}, and hence the inverse correlation length  $\xi_\beta$,
inherit the symmetries of $\bbZ^2$. 

\subsubsection{Multi-layer constrained macroscopic variational problems}
\label{subsub:CVP}
The relevant notions of (two-dimensional) Wulff shape and Wulff plaquettes are
introduced in Appendix~\ref{app:Wulff}. 
Given an area $\sfb\leq 4$, let $\xi_\beta (\sfb )$ be the minimal surface 
energy  over subsets of $[-1,1]^2$ of area $\sfb$.
In other words, $\xi_\beta(\sfb)$ is given by~\eqref{eq:Wulff-b} if $\sfb\in 
[0,w_\beta]$ and, accordingly, by~\eqref{eq:Wulff-P-b} if $\sfb\in [w_\beta,4]$.
Assuming that the expression for the contour weights 
in~\eqref{eq:Form-Weights} can be approximated by
$\sum_{i=1}^n \log q_\beta^{\sf f}(\gamma_i)$, that is, assuming that, at low 
temperatures, finite-volume effects and the interaction between different 
contours in the stack can be ignored, we conclude that information on the 
macroscopic scaling (by the linear system size $N$) of the stack of large 
contours should in principle be read from the constrained macroscopic 
variational problem

\smallskip 
\noindent
\textbf{(MCVP)}\qquad $\max_{\sfa_1\leq \sfa_2\leq\cdots\leq 4}
\bigl\{ -\sum_i \xi_\beta(\sfa_i) + \frac{1}{N} 
F_{\beta,N}(N^2\sfa_1,N^2\sfa_2,\ldots)\bigr\}$ .

\medskip 
Let us explore this variational problem for a class of examples we consider here. 
\smallskip 

\paragraph{\bf SOS model coupled to a  bulk Bernoulli field}
The multi-layer constrained variational problem \textbf{(MCVP)} was completely worked out 
in~\cite{IS16} for $F_{\beta,N }$ of the form~\eqref{eq:F-BSOS}, that is when  
\[ 
\frac{1}{N} F_{\beta, N}(N^2\sfa_1,N^2\sfa_2,\ldots)
=
-\frac{(\sfa - \sum_i \sfa_i)^2}{16D} . 
\]
In this case, the variation problem \textbf{(MCVP)} becomes 
\be{eq:VP-BSOS}
\min_{\sfa_1\leq\sfa_2\leq\ldots\leq 4}
\Bigr\{\sum_i \xi_\beta(\sfa_i) + \frac{(\sfa - \sum_i\sfa_i)^2}{16D}\Bigr\}.
\ee
Of course, there are different solutions for different values of $\sfa$. 
It is proven in~\cite{IS16} that~\eqref{eq:VP-BSOS} 
undergoes an infinite sequence of first-order phase transitions in the 
following sense: There exists a sequence $0=\sfb_0 <\sfb_1 <\sfb_2 <\ldots$, 
with $\lim_{k\to\infty}\sfb_k = \infty$, such that there is a unique solution 
to~\eqref{eq:VP-BSOS}, for any $\ell$ and for any $\sfa\in (\sfb_{\ell}, 
\sfb_{\ell +1})$. Furthermore, this solution contains exactly $\ell$ shapes  
\be{eq:SOS-bulk-stack} 
\sfS_1\subseteq \sfS_2 = \sfS_3 = \cdots = \sfS_\ell, 
\ee
and $\sfS_2,\ldots,\sfS_\ell$ are Wulff plaquettes, whereas $\sfS_1$ is 
either a Wulff plaquette identical to $\sfS_i$, $i\geq 2$, or it is a 
Wulff shape of the same radius as $\sfS_i$, $i\geq 2$. Moreover, starting 
from $\ell^* = \ell^*(\beta)$, the optimal solutions contain only Wulff 
plaquettes. 

\smallskip 

\paragraph{\bf Entropic repulsion}
Consider the expression~\eqref{eq:F-SOS} for $F_{N, \beta}$ (with the term
$f^*_0 \abs{\Lambda_N}$ dropped), applied to the microscopic areas
$\sfA_k = N^2 \sfa_k$: 
\be{eq:FN-entrop}
F_{N,\beta}(N^2\sfa_1,\ldots,N^2\sfa_n) 
= 
\bar{c}_\beta \sum_{k=1}^n {\rm e}^{-4\beta(n+1-k)} N^2\sfa_k
=
\bar{c}_\beta
N\sum_{k=1}^n {\rm e}^{-4\beta(n+1-k)} N\sfa_k.
\ee
The corresponding surface tension is $N\sum_{k=1}^n \xi_\beta (\sfa_k )$. 
This means that the top  layer (of area $N^2\sfa_1$) might appear only if 
\[
\max_{\sfa_1 \leq 4} \lbr \bar{c}_\beta {\rm e}^{-4\beta n}N\sfa_1 - \xi_\beta(\sfa_1)\rbr \geq 0 . 
\]
This falls into the framework of the dual constrained variational problem 
{\bf (DCVP)} discussed in Appendix~\ref{app:Wulff}.

Going back to~\eqref{eq:FN-entrop}, we conclude\footnote{Additional care is 
needed when $\nu_\beta$ is close to one of the end points of the above 
interval. We refer to~\cite{CLMST13} for precise statements and details.} 
that the number of layers $n^*$ should satisfy
\be{eq:SOS-nstar}
\nu_\beta \in \bigl( \bar{c}_\beta {\rm e}^{-4\beta(n^*+1)}N,
\bar{c}_\beta {\rm e}^{-4\beta n^*}N\bigr). 
\ee
Accordingly, define $\nu^1 = \bar{c}_\beta{\rm e}^{-4\beta n^*}N$ and
\be{eq:SOS-tilts} 
\nu^k
= 
{\rm e}^{4\beta(k-1)}\nu^1\text{, for $k=2, 3, \ldots$},
\ee
and consider the Wulff plaquettes
\be{eq:SOS-stack-var}
\sfP_\beta^{\sfb_1}\subset\sfP_\beta^{\sfb_2}\subset\cdots 
\text{ where }
\sfb_k = \sfa(\beta,\nu^k)
= 
\mathrm{argmax}_{\sfa\in [0,4]} \bigl\{ \nu^k\sfa - \xi_\beta(\sfa)\bigr\}. 
\ee

\smallskip
\paragraph{\bf Bulk magnetic field}
Consider~\eqref{eq:Form-Weights} and~\eqref{eq:F-SOS-B}. Recall that 
there are at most $\bar{n}$ contours in a stack. Restricting attention 
to contours with total length bounded above by  $K N$, we  infer that,
for all $\beta$ large enough,
\[ 
 \beta \sum_{i=1}^{\bar {n}} \abs{\gamma_i }
- 
\sum_{i=1}^{\bar{n}}\sum_{\calC} \Phi_{\beta,N}(\calC;\gamma_i)
- 
\Psi_{\beta,N}(\ugamma) \leq \bar{n} (K+1) N\beta. 
\]
Hence, 
\be{eq:bbeta-large}
\log q_{\beta,N}(\ugamma) \geq \sum_{\ell=1}^{\bar n} \lb r^*_{n-\ell +1}
N^2\sfa_\ell -  N (K+1)\beta \rb .
\ee
This means that if $r^*_{\bar n} >0$, the inner-most contour $\gamma_1$ tends to
fill in the whole box: $\sfa_1\to  4$, as $N\to\infty$. In other words, modulo a 
difficult  and, as we mentioned, still partially open question of characterization of 
$\bar{n} = \bar{n} (\beta,\lambda)$, 
the limiting variational problem in the case of a fixed bulk field $\lambda >0$ is 
somewhat trivial - all $\bar{n}$ limiting shapes are full $[-1, 1]^2$ squares. 

A more interesting situation should occur if the bulk field $\lambda = \lambda_N$ 
tends to zero as $N\to\infty$, but for the moment even a reliable  conjecture along these
lines seems to be beyond reach. 
\smallskip 

\paragraph{\bf Boundary magnetic fields.} 
The situation with the limiting macroscopic variational problem in the case
of boundary fields is somewhat similar to that in the bulk field. Namely, 
modulo an incomplete characterization of the typical, as $N\to\infty$, 
number of layers  $\bar{n}(\beta,h)$, the limiting variational problem
seems to be a trivial one. 
 
\subsubsection{Macroscopic scaling limits}
\label{subsub:msl}

All the results below are discussed under the tacit assumption that the 
inverse temperature $\beta$ is large enough.
 
Let us first consider SOS-type effective interface models with $\bbP_{N}^\beta$
as defined in~\eqref{eq:HamPM-SOS}. 

\smallskip 
\paragraph{\bf Entropic repulsion}

In the case of entropic repulsion, the scaling limits of large level lines 
were completely worked out in~\cite{CLMST13}. Let us formulate a particular 
instance of their results for sequences of side-lengths $N_\ell$ satisfying 
(recall~\eqref{eq:SOS-nstar} which determines $n^* = n^*(N)$)  
\be{eq:SOS-Nl-lim}
\lim_{\ell\to\infty} \bar{c}_\beta {\rm e}^{-4\beta n^*(N_\ell)} N_\ell
=
\nu^1>\nu_\beta,  
\ee
where $\nu_\beta$ is the critical value for the variational
problem~\textbf{(DCVP)}. 

\begin{theorem} 
\label{thm:SOS-lim} 
Assume~\eqref{eq:SOS-Nl-lim}. Let $\{\gamma_1 ,\gamma_2,\ldots\}$
be the ordered stack of large contours as in~\eqref{eq:contour-stack}. 
Then, using $\dd_{\sfH}$ for the Hausdorff distance,  
\be{eq:SOS-lim} 
\lim_{\ell\to\infty} \bbP_{N_\ell}^\beta \Bigl(
\dd_{\sfH} \bigl( \frac{\gamma_k}{N_\ell} , \partial\sfP_\beta^{\sfb_k} 
\bigr) > \epsilon \Bigr) = 0, 
\ee
for any $k\in\bbN$ fixed and any $\epsilon>0$, where $\sfb_k = 
\sfa(\beta,\nu^k)$ with $\nu^k$ defined as in~\eqref{eq:SOS-tilts}, and 
$\sfP_\beta^{\sfb_k}$ are optimal Wulff plaquettes as 
in~\eqref{eq:SOS-stack-var}. 
\end{theorem}

\smallskip
\paragraph{\bf SOS models with bulk and boundary fields} 
Formula~\eqref{eq:bbeta-large} indicates 
that the following should happen for $\beta$ large enough:  
\be{eq:SOS-lambda-lim}
\lim_{N\to\infty} 
\bbP_{N_\ell}^\beta \Bigl( \dd_{\sfH} 
\bigl( \frac{\gamma_1}{N} , \partial [-1, 1]^2 \bigr) > \epsilon \Bigr) = 0, 
\ee
for any $\epsilon >0$. 

Presumably, such results could be deduced from~\cite{DiMa94, CeMa96, 
LMazel96} for the range of $\beta$ and $\lambda$ to which the results of the 
latter papers apply.

The same is true regarding boundary fields in the regime in which the results
of~\cite{AlDuMi11} apply. 

More interesting phenomena might appear if one allows $\lambda_N\to\infty$ or 
$h_N\to\infty$ as the linear size $N$ of the system grows, but a rigorous
analysis will presumably require going well beyond the existing techniques. 

\smallskip 
\paragraph{\bf Facets in the SOS model coupled with bulk Bernoulli 
fields}

The scaling limits for large level lines under the measure $\bbQ_{N,\sfa}^\beta$ 
defined in~\eqref{eq:BSOS-Q} are studied in~\cite{IS16}. For the moment, the 
results there are formulated contingent to a proof of predominance of 
repulsion between different contours in a stack over weak interaction
between these contours, as described in more detail in 
Subsection~\ref{subsub:interaction} below. This issue was overlooked 
in~\cite{IS08}. Here is the \emph{conjectured} statement: 
 
\begin{theorem} 
There exists $\beta_0$ large enough such that the following holds: 
Fix $\beta >\beta_0$ and $\ell\in\bbN$. Let $\sfa\in (\sfb_{\ell},\sfb_{\ell+1})$,
where the sequence $\sfb_1,\sfb_2,\ldots$ is the one appearing in the solution of
the variational problem~\eqref{eq:VP-BSOS}. Then,
\be{eq:BSOS-claim1} 
\lim_{N\to\infty} 
\bbQ_{N,\sfa}^\beta(\text{Exactly one stack of $\ell$ large contours}) = 1 .
\ee
Furthermore, let $\sfS_1,\ldots,\sfS_\ell$ be the optimal shapes as 
in~\eqref{eq:SOS-bulk-stack}. Then, for any $\epsilon>0$ fixed, the unique 
$\ell$-stack $\{\gamma_1,\ldots,\gamma_\ell\}$ satisfies:
\be{eq:BSOS-claim2} 
\lim_{N\to\infty} 
\bbQ_{N,\sfa}^\beta \Bigl(
\max_{1\leq k\leq\ell} \dd_{\sfH} 
\bigl( \frac{\gamma_k}{N} , \partial\sfS_k \bigr) > \epsilon \Bigr) = 0. 
\ee
\end{theorem}
The results of~\cite{IofShlosTon15} justify the conclusion in the initial 
case of the first facet ($\sfb_1 < \sfa < \sfb_2$) as they imply that the 
interaction with the boundary of $\Lambda_N$ does not modify the surface 
tension. 
 
\subsection{Structure and fluctuations of interacting level lines}
\label{sec:EffectiveRandomWalk}

The macroscopic level lines are nested stacks of closed contours. A careful 
analysis of large-scale fluctuations of such stacks under the probability 
distribution with unnormalized weights~\eqref{eq:Form-Weights} should involve 
several steps: 
 
\noindent 
\step{1}
One should develop a fluctuation theory of mesoscopic segments of a single 
large contour. At this stage, the inverse correlation length  $\xi_\beta$ is incorporated 
and the effective local one-dimensional random walk structure of level lines 
is uncovered.
 
\noindent
\step{2}
One should develop a procedure, usually known as skeleton calculus, to patch 
mesoscopic segments into a single closed macroscopic contour.
 
\noindent
\step{3}
Different macroscopic contours in a nested stack are subject to entropic 
repulsion. One should check that the entropic repulsion prevails over the 
weak attraction due to the cluster sharing terms $\Psi_{\beta,N}$ 
in~\eqref{eq:Form-Weights}. In particular, the interaction between different 
contours should not modify the surface tension. Note that the multi-layer 
constrained variational problem \textbf{(MCVP)} is stated under the tacit
assumption that this indeed does not happen. 
 
\noindent 
\step{4}
Finally, one should understand what are the proper scaling regimes (as the 
size $N$ of the system goes to $\infty$) corresponding to the fluctuations
of the level lines under various area-type tilts $F_{\beta,N}$. 

\smallskip 
In the sequel we discuss what is known and what is not known along these 
lines. 
 
\subsubsection{Ising polymers and effective random walk representation} 
\label{sub:Ising} 

Let us discuss the low-temperature weights~\eqref{eq:polym-weights} for a 
mesoscopic segment $\gamma$ of a single level line.
This is a model of low-temperature Ising polymers and we shall closely follow 
the exposition in the recent work~\cite{IofShlosTon15}; in particular, we 
shall refer to Subsections~2 and~3 of the latter paper for missing details.
 
The low-temperature assumption comes into play through the very possibility 
to perform cluster expansions leading to~\eqref{eq:polym-weights}, and it is 
further quantified in terms of exponential decay properties of the cluster 
weights $\Phi_\beta(\calC;\gamma)$. Namely, assume that there exist some 
$\chi> 0$ such that, for all $\beta$ sufficiently large,
\be{eq:chi-decay}
\bigl| \Phi(\calC,\gamma) \bigr| < \exp\bigl\{ -\chi\beta
(\mathrm{diam} (\calC)+1) \bigr\}. 
\ee
Under~\eqref{eq:chi-decay}, the polymer weights $q^\sff_\beta(\gamma)$ 
in~\eqref{eq:polym-weights} can be rewritten in the following form 
(see~\cite[Subsection~3.2]{IofShlosTon15}): 
\be{eq:Wf-qf} 
q^\sff_\beta(\gamma)
=
\sum_{\underline{\calC}} {\rm e}^{-(\beta+c(\beta))\abs{\gamma}} 
\prod_i \rho_\beta(\calC_i,\gamma)
\df 
\sum_{\underline{\calC}} \rho^\sff_{{\beta}}
\bigl([ \gamma , \underline{\calC}]\bigr).
\ee
Above, the summation is with respect to all finite collections of connected 
clusters $\underline\calC = \{\calC_i\}$. The modified inverse temperature 
$\beta + c(\beta)$ is slightly larger than the original one:
$0 < c(\beta) < {\rm e}^{-\chi\beta}$. Finally, the product cluster weights 
$\rho_\beta(\calC,\gamma)$ are \emph{non-negative} and exponentially decaying:
\be{eq:rho-weights} 
0\leq \rho_\beta(\calC,\gamma) \leq {\rm e}^{-\chi
\beta
(\mathrm{diam} (\calC)+1)}\IF{\text{$\calC$ is incompatible with $\gamma$}}. 
\ee
The weights $\rho_\beta(\calC,\gamma)$ are translation invariant. 
In such circumstances, a straightforward adjustment of sub-additivity 
arguments implies that the limit in~\eqref{eq:tau-beta} exists for all 
$\beta$ sufficiently large, and that the corresponding surface tension 
$\xi_\beta$ is strictly positive and convex.

Furthermore, much sharper results hold: Let us call couples 
$\fra = [\gamma,\underline{\calC}]$ \emph{animals}. Given a cone $\calY$ and 
a point $\sfx\in\gamma = (\sfx_0,\ldots,\sfx_n)$, let us say that $\sfx$ is 
a $\calY$-break point of $[\gamma,\underline{\calC}]$ if it is an interior 
point of $\gamma$ (that is, $\sfx\in\{\sfx_1,\ldots.\sfx_{n-1}\}$) 
and the following holds:
\be{eq:y-bp} 
\{\sfx_0,\ldots,\sfx_i\} \subset (\sfx - \calY),
\quad
\{\sfx_i,\ldots, \sfx_n\} \subset (\sfx + \calY)
\text{ and }
\cup_i \calC_i \subset (\sfx - \calY)\cup(\sfx + \calY).
\ee
Evidently, if $\sfx$ is a $\calY$-break point of an animal 
$[\gamma,\underline{\calC}]$, then we can represent it as a concatenation
\[ 
\fra
= [\gamma,\underline{\calC}]
= [\gamma^1,\underline{\calC}^1]\circ [\gamma^2,\underline{\calC}^2]
= \fra^1\circ\fra^2
\]
such that $\gamma^1,\underline{\calC}^1\subset (\sfx - \calY)$ and 
$\gamma^2,\underline{\calC}^2\subset (\sfx +\calY)$. 
\begin{definition} 
\label{def:irreducible}
Given a cone $\calY$, let us say that an animal $\fra = 
[\gamma,\underline{\calC}]$ with $\gamma = (\sfx_0,\sfx_1,\ldots,\sfx_n)$ is 
$\calY$-irreducible if it does not have break points.
\end{definition}
Unnormalized Wulff shapes $\calW_\beta$ are defined 
in~\eqref{eq:def-Wulff} of Appendix~\ref{app:Wulff}. 
With each $\sfh\in\partial\calW_\beta$, we can 
associate a convex cone
\be{eq:Y-h} 
\calY_\sfh = \setof{\sfx}{\sfh\cdot\sfx \geq \epsilon}.
\ee
We assume that $\epsilon > 0$ is sufficiently small, so that $\calY_\sfh$
always contains a lattice direction in its interior (and hence there are 
paths $\gamma$ which satisfy $\gamma\subset \calY_\sfh$). 

\smallskip
\paragraph{\bf Ornstein--Zernike theory}
The relevant input from the OZ theory (see for instance 
\cite[Subsections~3.3 and~3.4]{IV08} and~\cite[Section~4.1]{IofShlosTon15}) 
can be summarized as follows: For $\sfh\in\partial\calW_\beta$, let 
$\calY_\sfh$ be the cone defined in~\eqref{eq:Y-h} and let $\sfA_\sfh$
be the corresponding set of irreducible animals. For $\fra = 
[\gamma,\underline{\calC}]\in\sfA_\sfh$ with $\gamma = 
(\sfx_0,\ldots,\sfx_n)$, define the length $\abs{\sfa} = \abs{\gamma}$ and 
the displacement $\sfX(\sfa) = \sfx_n-\sfx_0$, and set 
(recall~\eqref{eq:W-beta-1} and~\eqref{eq:Wf-qf})
\be{eq:OZ-weights}
\bbP_\beta^{\sfh }(\fra)
= 
{\rm e}^{ \sfh\cdot\sfX(\sfa)} \rho_\beta^\sff(\fra).
\ee
Let us say that an $\calY_\sfh$-irreducible animal is a left diamond, 
respectively a right diamond, if 
\be{eq:irreducible} 
\gamma \subset (\sfx_n - \calY),
\text{ respectively }
\gamma \subset (\sfx_0 + \calY).
\ee
An animal $\fra$ is said to be a diamond if it is both a left and a right 
diamond. We use the notation $\sfD^\sfl_\sfh, \sfD^\sfr_\sfh$ and $\sfD_\sfh 
= \sfD^l_\sfh\cap\sfD^r_\sfh$ for the corresponding sets of animals. 
\begin{theorem}
\label{thm:OZ}
For all $\beta$ large enough and for every $\sfh \in \partial\calW_\beta$, 
$\bbP_\beta^{\sfh }$ is a probability distribution on $\sfD_\sfh$ with 
exponentially decaying tails on $\sfA_\sfh$. That is,  
\be{eq:exp-tails} 
\sum_{\fra\in\sfD_\sfh} \bbP_\beta^{\sfh}(\fra) = 1
\text{ and }
\sum_{\fra\in\sfA_\sfh} \bigl\{ 
\bbP_\beta^{\sfh}(\abs{\fra} > k)
+ 
\bbP_\beta^{\sfh}(\sfX(\fra) > k) \bigr\}
< {\rm e}^{-\nu_\beta k} ,
\ee
for any $k>1$, where $\nu_\beta$ does not depend on $\sfh$.
 
Furthermore, $\partial \calW_\beta$ is locally analytic with uniformly 
strictly positive curvature . In fact, the parametrization of $\partial 
\calW_\beta$ in a neighborhood of $\sfh$ can be described as 
\be{eq:p-Wbeta-eq} 
\sfg +\sfh \in \partial\calW_\beta
\quad\Leftrightarrow\quad 
\bbE_\beta^\sfh \bigl( {\rm e}^{ \sfg\cdot\sfX(\fra)}\bigr) = 1.
\ee
In particular, $\xi_\beta$ is differentiable at any $\sfx\neq 0$ and 
\be{eq:sfh-sfx}
\sfh_\sfx =  \nabla\xi_\beta(\sfx)
\ee
is the unique point on $\partial\calW_\beta$ satisfying 
$\sfh_\sfx\cdot\sfx = \max_{\sfh\in\partial\sfK_\beta} \sfh\cdot\sfx = 
\xi_\beta(\sfx)$. 
\end{theorem} 
Theorem~\ref{thm:OZ} indicates that a typical mesoscopic segment of a level 
line should look like a one-dimensional necklace of irreducible animals. Let 
us make this precise. 

\smallskip 
\paragraph{\bf Effective RW representation of mesoscopic segments of 
level lines}

Let $\sfx$ be a distant point and consider the set of all animals $\fra = 
[\gamma,\underline{\calC}]$ from the origin to $\sfx$. In other words, 
consider the set of all animals $\sfa$ with $\sfX(\fra)=\sfx$. 

Now, \eqref{eq:tau-beta} implies that $q^\sff_\beta (\sfx ) \asymp 
{\rm e}^{-\xi_\beta(\sfx)}$. Consider $\sfh = \sfh_\sfx$ as defined 
in~\eqref{eq:sfh-sfx}. Then,
\[
1 \asymp {\rm e}^{\xi_\beta(\sfx)} q^\sff_\beta(\sfx)
= 
\sum_{\sfX(\sfa)=\sfx} {\rm e}^{\sfh\cdot\sfX(\sfa)} 
\rho_\beta^\sff(\fra).
\]
By~\eqref{eq:OZ-weights} and~\eqref{eq:exp-tails}, and up to 
corrections of order ${\rm e}^{-\nu_\beta \|\sfx\|}$, we can restrict 
attention to animals $\fra$ of the form 
\be{eq:a-concat} 
\fra = \frb^\sfl\circ\fra^1\circ\cdots\circ\fra^m\circ\frb^\sfr, 
\ee
with $\frb^l\in\sfD^l_\sfx$, $\frb^r\in\sfD^r_\sfh$ and $\sfa^i\in\sfD_\sfh$. 

\begin{figure}[t]
\begin{center}
\includegraphics[width=9cm]{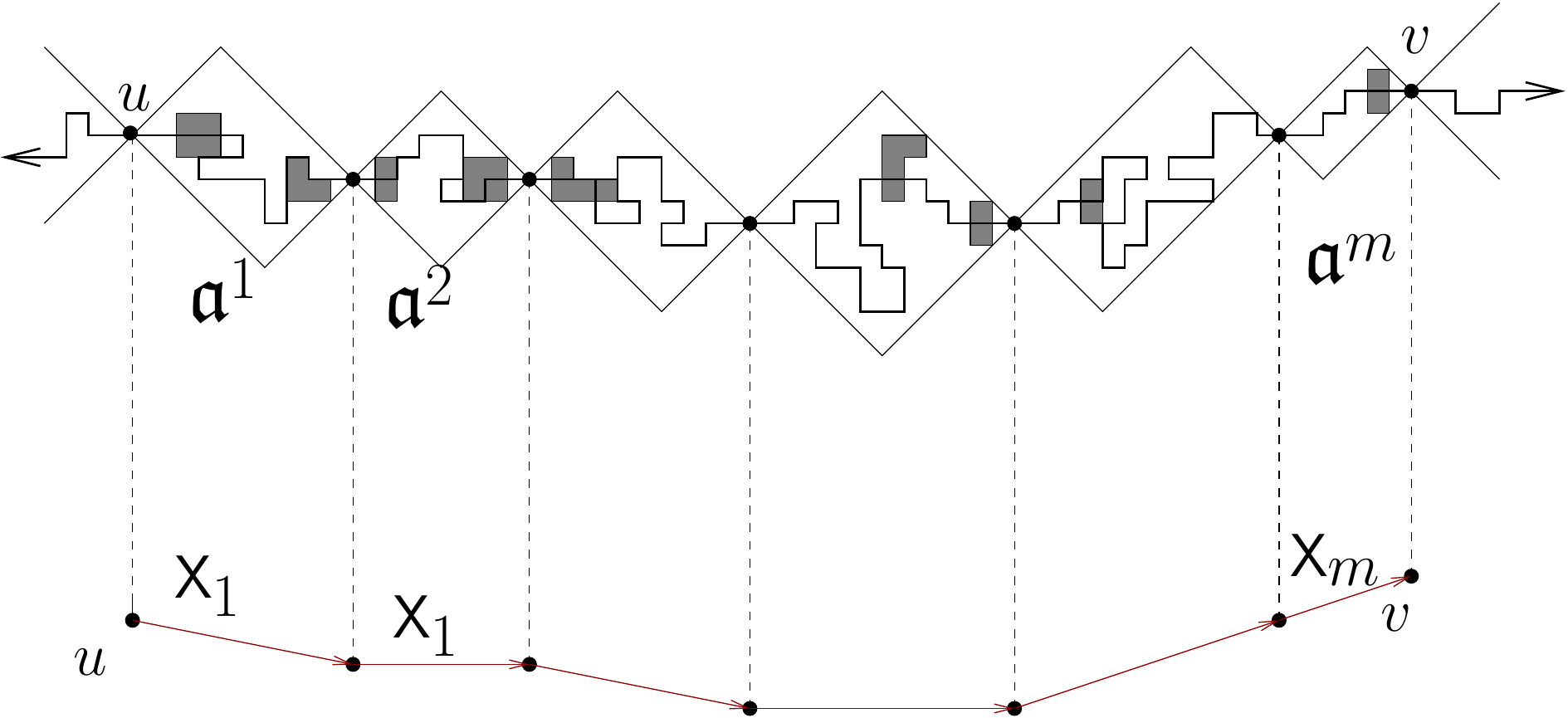}
\end{center}
\caption{Effective random walk structure of connection between $u$ and $v$.}
\label{fig:RWAnimal}
\end{figure}

The effective random walk representation, see Figure~\ref{fig:RWAnimal}, 
of mesoscopic segments 
$\gamma : 0\mapsto \sfx$ can be recorded as follows: Let 
$\sfX_1,\sfX_2,\ldots$ be i.i.d.\ $\calY_\sfh$-valued random variables 
distributed according to $\bbP_\beta^{\sf}$, that is,
\be{eq:X-distr}
\bbP_\beta^{\sfh}(\sfX_i  = \sfx)
=
\sumtwo{\fra\in\sfD_\sfh}{\sfX(\fra)=\sfx} \bbP_\beta^{\sfh}(\fra).
\ee
Set $\sfS_j = \sum_{i=1}^j \sfX_i$. Then $\gamma$ is approximated by the 
polygonal line, or equivalently, the trajectory of the effective random walk  
through the vertices 
\[
0,\, \sfX(\frb^\sfl),\,  \sfX(\frb^\sfl) + \sfX_1, \ldots,\,
\sfX(\frb^\sfl) + \sfS_m,\, \sfx .
\]

In this way, the probability distribution of the effective random walk comes 
from the normalization of the partition function 
\be{eq:sum-segments}
{\rm e}^{\xi_\beta(\sfx)}
q^\sff_\beta (\sfx) \bigl( 1 + {\rm o}({\rm e}^{-\nu_\beta \|\sfx\|})\bigr)
= 
\sumtwo{\frb^\sfl\in\sfD^\sfl_\sfh}{\frb^\sfr\in\sfD^\sfr_\sfh}
\bbP_\beta^{\sfh}(\frb^\sfl) \bbP_\beta^{\sfh}(\frb^\sfr)
\sum_{m=1}^\infty \otimes \bbP_\beta^{\sfh}\bigl(
\sfX(\frb^\sfl) +  \sfS_m +\sfX(\frb^\sfr) = \sfx \bigr).
\ee
The collection of clusters $\underline{\calC}$ in $\fra = 
[\gamma,\underline{\calC}]$ should  be viewed as \emph{hidden variables}. The 
real object of interest are paths  $\gamma \df \gamma(\fra)$. 
For $\fra\in\sfD^\sfl_\sfh,\sfD^\sfr_\sfh,\sfD_\sfh$, such paths $\gamma = 
(\sfx_0,\ldots,\sfx_n)$ belong, respectively, to the sets $\calP^\sfl_\sfh$, 
$\calP^\sfr_\sfh$ and $\calP_\sfh = \calP^\sfl_\sfh\cap\calP^\sfr$, where 
\be{eq:P-sets} 
\calP^\sfl_\sfh = \setof{\gamma}{\gamma\subset \sfx_n - \calY_\sfh}
\text{ and } 
\calP^\sfr_\sfh = \setof{\gamma}{\gamma\subset \sfx_0 + \calY_\sfh}.
\ee
Accordingly, \eqref{eq:OZ-weights} gives rise to a probability distribution 
on $\calP_\sfh$ and to finite measures with exponentially decaying tails on 
$\calP^\sfl_\sfh$ and $\calP^\sfr_\sfh$, which we, with a slight abuse of 
notation, continue to denote $\bbP_\beta^{\sfh}$. For instance, for 
$\gamma\in\calP_\sfh$, 
\be{eq:Prob-h-gamma} 
\bbP_\beta^{\sfh }(\gamma) = \sum_{\fra \,:\, \gamma(\fra)=\gamma} 
\bbP_\beta^{\sfh}(\fra) \IF{\fra\in\sfD_\sfh}.
\ee
In this way, instead of~\eqref{eq:a-concat}, we can consider 
\be{eq:a-concat-gamma} 
\gamma = \eta^\sfl\circ \gamma^1\circ\dots \circ \gamma^m \circ \eta^\sfr
\ee
and, accordingly, instead of~\eqref{eq:sum-segments}, we can write 
\be{eq:sum-segments-gamma}
{\rm e}^{\xi_\beta(\sfx)}
q^\sff_\beta(\sfx) \bigl( 1 +{\rm o}({\rm e}^{-\nu_\beta \|\sfx\|})\bigr)
= 
\sumtwo{\eta^\sfl\in\calP^\sfl_\sfh}{\eta^\sfr\in\calP^\sfr_\sfh}
\bbP_\beta^{\sfh}(\eta^\sfl) \bbP_\beta^{\sfh}(\eta^\sfr)
\sum_{m=1}^\infty \otimes \bbP_\beta^{\sfh} \bigl( \sfX(\eta^\sfl) + \sfS_m 
+ \sfX(\eta^\sfr) = \sfx \bigr).
\ee

\subsubsection{Skeletons and interaction between different contours}
\label{subsub:interaction}

The formula~\eqref{eq:sum-segments-gamma} furnishes a probabilistic 
description of an open ``linear'' portion $\gamma$ of a level line between 
two distant points $\sfx$ and $\sfy$. 
When $\|\sfy-\sfx\|$ is large, it both recovers the macroscopic inverse correlation
length 
$\xi_\beta(\sfy-\sfx)$ and indicates the fluctuation structure of 
$\gamma$ in the corresponding reduced ensembles. Recall, however, 
that the microscopic level lines are closed contours of size $N$. The idea of 
skeleton calculus is to go to an intermediate coarse-graining scale, say 
$N^\alpha$ for $\alpha\in (0,1)$, and try to study closed contours of size 
$N$ as a concatenation of open paths of size $N^\alpha$. This is with a hope 
that patching such open paths together will lead to controllable corrections 
to asymptotic formulas such as~\eqref{eq:sum-segments-gamma}. In the latter 
case, since the number of different $N^\alpha$-skeletons is bounded above by  
${\rm e}^{cN^{1-\alpha}}$, one infers concentration near macroscopic shapes 
of minimal $\xi_\beta$-surface energy. We refer to the groundbreaking 
book~\cite{DKS} where this idea was introduced in the context of the 
low-temperature 2D Ising model, and to various implementations of skeleton 
calculus in subsequent works~\cite{ACC90, P91, PfVe1997, ISc98, BoIoVe2000, 
Al01, BCK03, IS08}.

The implementation of the above procedure should hinge on an argument which 
would imply that different mesoscopic segments do not interact, in the sense 
that the total surface energy is the undistorted sum of surface energies over 
different mesoscopic segments. The same applies if one considers several 
macroscopic level lines.

Let us elaborate on this problem. Let 
$\{\sfx_1,\sfy_1\},\ldots,\{\sfx_k,\sfy_k\}$ be a collection of pairs of 
points in $\bbZ^2$, which represent neighboring  vertices of skeletons of the 
same or different level lines. 
That is, we assume that $\|\sfy_j-\sfx_j\|\gg 1$.
To fix ideas and in order to avoid redundant notation, let us assume that all 
$\sfx_i$-s lie on a vertical axis through $-\ell\sfe_1$, that is, 
$\sfx_i\cdot\sfe_1 \equiv -\ell$, and similarly that $\sfy_i$-s lie on a 
vertical axis through $\ell\sfe_1$. Assume also that these vertices are 
ordered in the sense that
\[
\sfx_1\cdot\sfe_2 \leq \sfx_2\cdot\sfe_2 \leq \cdots \leq \sfx_k\cdot \sfe_2, 
 \]
and the same regarding $\sfy_i$-s. Let $\underline{\gamma} = 
(\gamma_1,\ldots,\gamma_k)$ be a family of compatible paths $\gamma_i :\sfx_i 
\to \sfy_i$. Recall that the notion of compatibility comes from the splitting 
rules employed in the construction of ordered stacks of microscopic level lines 
as indicated in Subsection~\ref{susub:level-lines}. In particular, the
$\gamma_i$-s do not cross. According to~\eqref{eq:Form-Weights} (for the 
moment, we ignore the confining geometry of $\Lambda_N$ as well as the area 
tilt $F_{\beta,N}$), the free joint weights $\log q_{\beta}^\sff(\ugamma)$
are given by 
\be{eq:Form-Weights-k}
\log q_{\beta}^\sff(\ugamma)
=
- \beta \sum \abs{\gamma_i}
+ \sum_{i=1}^k\sum_{\calC} \Phi_{\beta}(\calC;\gamma_i)
+ \Psi_{\beta} (\ugamma)
=
\sum_1^k \log q_{\beta}^\sff(\gamma_i)  + \Psi_{\beta}(\ugamma),
\ee
where the interaction term $\Psi_{\beta}$ is due to overcounting (cluster 
sharing) and is given by
\be{eq:Psi-beta} 
\Psi_{\beta}(\ugamma)
=
\sum_{\calC} \bigl\{
\Phi_{\beta}(\calC,\ugamma) - \sum_1^k \Phi_{\beta}(\calC;\gamma_i)
\bigr\}.
\ee
The question is whether one can control the partition functions
\be{eq:PF-uxuy} 
q_{\beta}^\sff(\ul{\sfx},\ul{\sfy})
\df
\sum_{\ugamma\text{ compatible}} q_{\beta}^\sff(\ugamma)
=
\sum_{\ugamma\text{ compatible}} {\rm e}^{\Psi_{\beta}(\ugamma)} 
\prod_{i=1}^k q_{\beta}^\sff(\gamma_i)  
\ee
in terms of $\prod_{i=1}^k q_{\beta}^\sff(\sfy_i -\sfx_i) \asymp 
{\rm e}^{-\sum_i \xi_\beta(\sfy_i -\sfx_i)}$.

Coarse  lower bounds are easy: 
Indeed, by~\eqref{eq:sum-segments-gamma}, one can confine paths $\gamma_i$ to
distant tubes, and in the latter case the exponential decay 
in~\eqref{eq:chi-decay} renders negligible the contribution due to the 
interaction.

Upper bounds are more difficult. 
In view of the effective random walk 
representation~\eqref{eq:sum-segments-gamma},
the following refined version of the above question makes sense. Consider 
the partition functions
\be{eq:PF-uxuy-not} 
\widehat{q}^{\;\sff}_\beta (\ul{\sfx},\ul{\sfy})
\df
\sum_{\ugamma\text{ compatible}} \prod_{i=1}^k q_{\beta}^\sff(\gamma_i)
\ee
and compare the distribution of $\ugamma$ under $q_{\beta}^\sff$ with 
the distribution of $\ugamma$ under $\widehat{q}^{\;\sff}_\beta$. 

The point is that, under $\hat{q}^{\;\sff}_\beta$, we are essentially 
talking about $k$ ordered effective random walks, and the behavior of the 
latter, including entropic repulsion, is essentially well understood. 

The competition between a potential attraction through the $\Psi_\beta$ term 
in~\eqref{eq:PF-uxuy} and the entropic repulsion between the paths $\gamma_i$ 
is highly non-trivial. It is true that, under~\eqref{eq:chi-decay}, the 
interaction $\Psi_\beta$ has an exponential decay of order $\beta$, but so is 
the variance of effective random walk steps in~\eqref{eq:X-distr}, at least 
in lattice directions. Furthermore, an additional complication comes with 
attempts to derive bounds which would hold uniformly up to leading terms in 
the number of paths $k$. 

In fact~\cite{IofShlosTon15}, in the full generality of Ising polymers, just 
having $\chi>0$ in~\eqref{eq:chi-decay} is not enough to ensure that entropic 
repulsion wins. The threshold value is conjectured to be $\frac{1}{2}$. 
In~\cite{IofShlosTon15}, it is shown that, for $\chi > \frac{1}{2}$, a half 
space polymer which interacts with a hard wall eventually behaves as the free 
full-space polymer subject to entropic repulsion. An ad hoc counter-example 
for $\chi =\frac{1}{2}$ is also constructed. A solution to the general 
question about $k$ interacting paths as stated above is, for the moment, not 
written down. This issue will be addressed in  the forthcoming 
work~\cite{IS16}. 

In the case of entropic repulsion~\cite{CMT14,CLMST13}, upper bounds follow 
from specific FKG properties of the underlying SOS-Hamiltonians. Namely, for 
 (a slight finite volume modification of) 
the partition functions $q_{\beta}^\sff(\ul{\sfx},\ul{\sfy})$ defined 
in~\eqref{eq:PF-uxuy}, the following holds:
\be{eq:PF-upper-SOS} 
q_{\beta}^\sff (\ul{\sfx},\ul{\sfy})
\leq 
\prod_{i=1}^k q_{\beta}^\sff(\sfy_i -\sfx_i).
\ee

\subsubsection{Scaling and fluctuations under area tilts} 

In all the examples discussed in Subsection~\ref{subsub:msl}, the limiting 
shapes for large level lines are either Wulff shapes or Wulff plaquettes. 
Actually, in the models we consider, Wulff shapes appear only as top droplets
in the case of an SOS surface coupled with Bernoulli  bulk fields, and only when 
the number of facets is bounded by some $\ell^*(\beta)$; see the discussion 
just after~\eqref{eq:SOS-bulk-stack}. In the sequel, we shall consider the
fluctuations of large level lines around the limiting {\em plaquettes}. 
The latter contain macroscopic flat segments on the boundary $\partial [-1,1]^2$,
and we shall focus on the fluctuations of microscopic level lines 
$\gamma_1,\gamma_2,\ldots$ away from these segments. 

To fix ideas, let $[-\delta,\delta]\times\{-1\}$ be a flat segment of the 
top limiting plaquette which appears in either of the scaling limits 
discussed in Subsection~\ref{subsub:msl}. Then it appears as a flat segment 
in all the subsequent limiting shapes. Let us zoom in near the microscopic 
boundary of $\partial\Lambda_N$ which contains this segment. We want to study 
the behavior of microscopic portions $\eta_1,\eta_2,\ldots,\eta_n$ of the 
level lines $\gamma_1,\gamma_2,\ldots,\gamma_n$ above this segment. 
For the sake of this discussion, we shall consider a somewhat simplified 
picture along the lines of Subsection~\ref{subsub:interaction}.

\begin{figure}[t]
\begin{center}
\includegraphics[width=8.7cm]{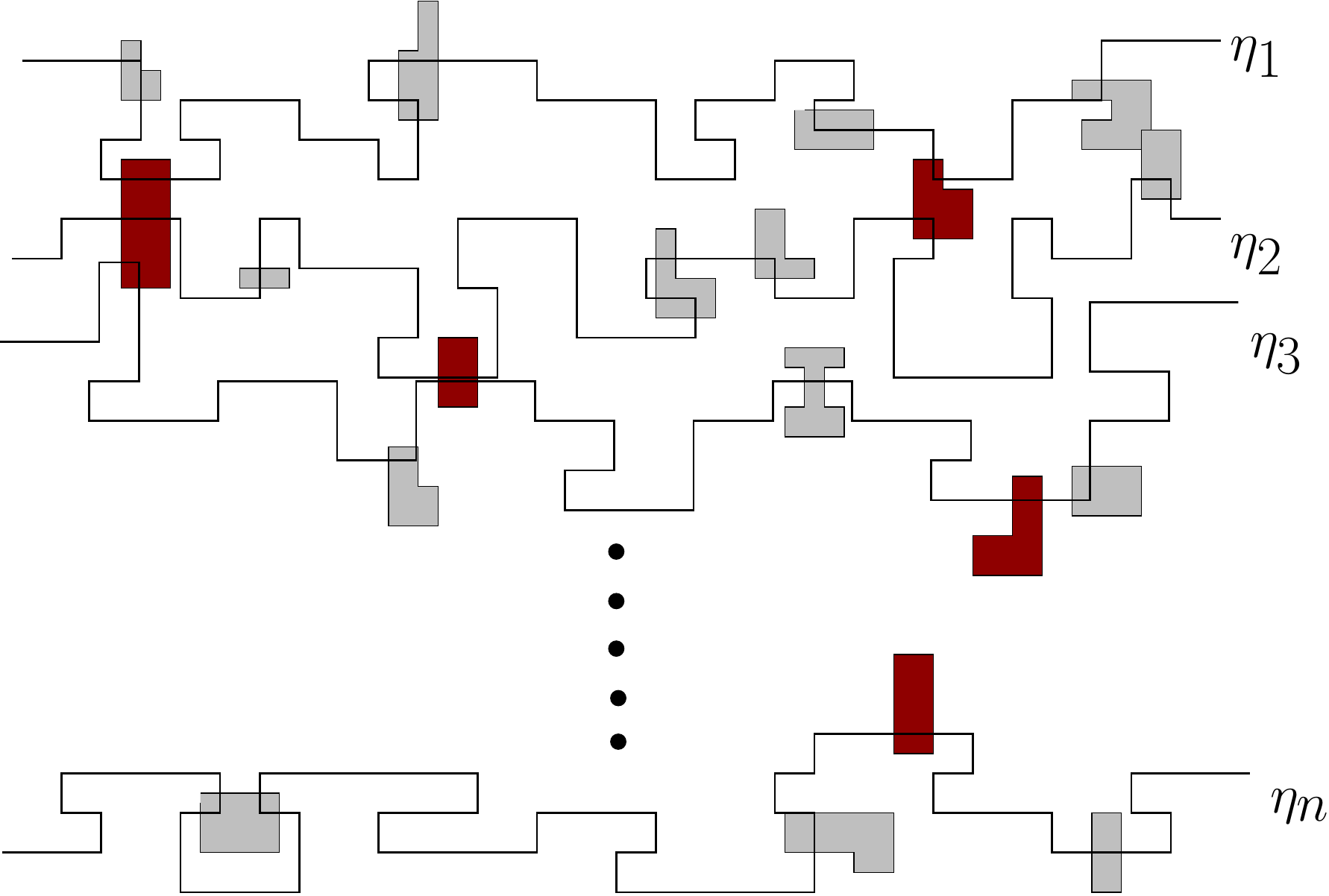}
\end{center}
\caption{Interaction by cluster sharing between pieces of level lines 
$\eta_n, \dots , \eta_1$.}
\label{fig:Interaction}
\end{figure}

Let $\{\sfx_1,\sfy_1\}, \ldots, \{\sfx_n,\sfy_n\}$ be a collection of pairs 
of points in the upper half lattice $\PHS^2 = \setof{\sfx\in 
\bbZ^2}{\sfx\cdot\sfe_2 \geq 0}$ such that $\sfx_i\cdot\sfe_1 \equiv -\delta 
N$, $\sfy_i\cdot\sfe_1 \equiv \delta N$ and
\be{eq:x-order}
\sfx_1\cdot\sfe_2 \geq \sfx_{2}\cdot\sfe_2 \geq \cdots \geq \sfx_n\cdot\sfe_2 
\geq 0 ,
\ee
with the same vertical ordering holding for  $\sfy_1,\ldots,\sfy_n$. Thus, 
$\underline{\sfx}$ and $\underline{\sfy}$ represent collections of initial 
and final vertices of \emph{ordered} portions $\eta_1,\ldots,\eta_n$ of the
level lines $\gamma_1,\ldots,\gamma_n$. Let us ignore the interactions 
between the $\eta_i$-s and the remaining pieces of level lines 
$\gamma_1\setminus\eta_1,\ldots,\gamma_n\setminus\eta_n$. That is, we 
assume that the \emph{free} weights of the paths $\eta_i$-s are given by a 
modification of~\eqref{eq:Form-Weights-k}:
\be{eq:Form-Weights-k-plus}
\log q_{\beta}^+(\ueta)
=
- \beta \sum \abs{\eta_i }
+ \sum_{i=1}^k\sum_{\calC\subset\PHS^2} \Phi_{\beta}(\calC;\eta_i)
+ \Psi_{\beta}^+(\ueta)
=
\sum_{i=1}^k \log q_{\beta}^+(\eta_i) + \Psi_{\beta}^+(\ueta), 
\ee
where the $+$-superscript indicates that we are taking into account the 
interaction with the boundary of $\Lambda_N$. The term $\Psi_{\beta}^+(\ueta)$
represents multi-body interactions by cluster sharing between $\eta_i$-s, see 
Figure~\ref{fig:Interaction}. 

The main issue here is to understand what type of corrections one should add 
to the effective ensemble above in order to incorporate fluctuations away 
from the rescaled macroscopic limiting shapes. In the sequel, we shall assume 
that all $\eta_k$-s are confined to the strip $[-\delta N , \delta N]\times 
\bbZ_+$ and use $\sfB_k$ to denote the area below $\eta_k$. 

\smallskip 
\paragraph{\bf SOS model coupled with Bernoulli bulk fields}

By assumption, the rescaled optimal limiting plaquettes
$N\sfP^{\sfa_1}_\beta = N\sfP^{\sfa_2}_\beta = \cdots = N\sfP^{\sfa_n}_\beta$ 
of areas $\sfA_1 = N^2 \sfa_1, \ldots, \sfA_n = N^2 \sfa_n$ stick to the 
boundary of $\partial \Lambda_N$ along the segment in question. So, in the 
rescaled picture of optimal shapes, facets pile up along flat segments of the 
boundary and $\sum_{k=1}^n \sfA_k$ is the optimal
total volume below the three-dimensional surface. The level lines 
$\eta_1,\ldots,\eta_n$ fluctuate away from this flat boundary and, as a 
result, reduce this optimal volume by $\sum_{k=1}^n \sfB_k$. This deficit 
must be compensated by excess fluctuations of the bulk Bernoulli field. 
Consider now~\eqref{eq:F-BSOS}. The excess price is, up to lower-order terms, 
given by
\be{eq:BSOS-excess} 
- \frac{(\sfa N^2 -\sum_i \sfA_i)^2}{2D \abs{B_N}}
+ \frac{(\sfa N^2 -\sum_i \sfA_i + \sum_i\sfB_i)^2}{2D \abs{B_N}}
\approx \frac{c_\beta(\sfa)}{N} 
\sum_{i=1}^n \sfB_i.
\ee
Therefore, in the case of an SOS model coupled with $3d$ Bernoulli bulk 
fields, the large-scale behavior of $n$ ordered segments $\eta_1, \ldots, 
\eta_n$ of large macroscopic level lines should be captured in the context of 
the following effective model: 
Fix $\delta >0$ and $c>0$; consider the weights
\be{eq:BSOS-n-weights}
\log q_{\beta,\delta,c,N}^{+,n}(\ueta)
=
\sum_{i=1}^k \log q_{\beta}^+(\eta_i)
+ \Psi_{\beta}^+ (\ueta) - \frac{c}{N}\sum_{i=1}^n \sfB_i
\ee
and let $\bbQ_{\beta,\delta,c,N}^{+,n}$ denote the corresponding probability 
distribution. 

\begin{conjecture}
\label{conj:BSOS}
For any $n$ fixed, the paths $\eta_i$ should be rescaled by $N^{2/3}$ in the 
horizontal direction and by $N^{1/3}$ in the vertical direction. Under this 
rescaling, the limit of $\bbQ_{\beta,\delta,c,N}^{+,n}$ does not depend on 
$\delta$ and is given in terms of ergodic Dyson Ferrari--Spohn diffusions, 
which we shall describe in  Appendix~\ref{app:FS}. 
\end{conjecture}

\smallskip
\paragraph{\bf Entropic repulsion}

Let $N_\ell$ be a diverging sequence of linear sizes satisfying 
assumption~\eqref{eq:SOS-Nl-lim}. Theorem~\ref{thm:SOS-lim} implies that 
the rescaled level lines should concentrate around the optimal plaquettes as 
in~\eqref{eq:SOS-stack-var}. Again, these optimal shapes pile up near the 
flat pieces along the boundary of $\Lambda_{N_\ell}$. Let us find out the 
appropriate expressions for the extra cost associated to the fluctuations of 
$\eta_1,\ldots,\eta_{n^*}$ away from this boundary. At the level of the 
variational problem, the surface is at height $n^*$. However, 
microscopically, the area of the surface at height $(n^*-k )$ is $(\sfB_{k} - 
\sfB_{k +1})$. In view of~\eqref{eq:SOS-lop}, this entails (as usual, up to 
lower-order corrections) an extra price
\begin{multline}
\label{eq:extra-SOS}
c_\beta \sum_{k=1}^{n^*-1} (\sfB_{k} - \sfB_{k+1}) \lb {\rm e}^{-4\beta(n^*-k+1)} 
- {\rm e}^{-4\beta(n^*+1)}\rb 
+ c_\beta \sfB_{n^*} \lb {\rm e}^{-4\beta} - {\rm e}^{-4\beta(n^*+1)}\rb 
\\
= 
c_\beta \lb 1 - {\rm e}^{-4\beta}\rb = 
\bar{c}_\beta \sum_1^{n^*} {\rm e}^{-4\beta(n^*-k+1)} \sfB_k  
\stackrel{\eqref{eq:SOS-Nl-lim}}{=} 
\sum_1^{n^*} \frac{\nu^1 (1+\smo{1} ){\rm e}^{4\beta(k -1)}}{N_\ell} \sfB_k.
\end{multline}
In other words, in the case of entropic repulsion, not only does the number 
of paths grow logarithmically with the linear size $N$, but also the tilt is 
growing exponentially with the position $k$ of $\eta_k$ in the ordered stack. 

In~\cite{CLMST13}, it is proven that fluctuations of the top segment are 
bounded above by $N^{1/3 +\epsilon}$ for any $\epsilon>0$.
At this stage, it is not clear whether $N^{1/3}$ is the correct scaling of 
the paths $\eta_k$ in vertical direction, or whether additional, for instance 
logarithmic, corrections are needed. Furthermore, there are no clear 
conjectures about the existence and the structure of fluctuation scaling 
limits for the whole stack. 

\smallskip 
\paragraph{\bf Bulk and boundary fields}
Let $\lambda>0$ be the bulk field. According to the discussion in 
Subsection~\ref{subsub:msl}, one should expect a stack of roughly $\bar{n} (\beta , \lambda)$ 
large 
macroscopic level lines, where $\bar{n}$ is given by~\eqref{eq:SOS-heights-n}. As 
$N$ goes to infinity, these contours should stick to the boundary 
$\partial\Lambda_N$ with a presumably bounded order of fluctuations, 
depending on $\lambda$ but not on $N$.
An interesting phenomenon, which bears a resemblance to what happens in 
the case of entropic repulsion, should appear if we let the bulk 
filed vanish; 
$\lambda  = \lambda_N\searrow 0$,  as the linear size of the system $N$. 
As we have already mentioned, presumably one has to go well beyond existing 
techniques in order to study this issue. The same regarding boundary fields.

\section{Non-colliding random walks under generalized area tilts}
\label{sec:RW}

In this Section we shall explain the rationale behind Conjecture~\ref{conj:BSOS}, 
which is based on considering a simplified random walk model for ~\eqref{eq:BSOS-n-weights}.
The simplification will be two-fold: First of all, we shall
ignore the interaction term $\Psi_\beta^+ (\ueta)$ and keep only the hard-core
constraint: paths $\eta_i$ will still be ordered. Furthermore, irreducible pieces
in~\eqref{eq:sum-segments-gamma} will have unit horizontal span or, in other 
words, will be represented by independent steps of a one-dimensional random walk.
On the other hand, we shall consider self-potentials of a more general form than the
linear area tilts appearing in the last term on the right hand side of~\eqref{eq:BSOS-n-weights}. 

Let $\lbr \sfx_1 , \sfy_1\rbr, \dots , \lbr \sfx_n , \sfy_n\rbr$ be a collection 
of ordered vertices  satisfying \eqref{eq:x-order}. To simplify the exposition  we 
shall assume that $\delta = 1$, that is we shall assume that for any 
$i=1, \dots , n$, the scalar products $\sfx_i\cdot\sfe_1 =-N$ and $\sfy_i\cdot\sfe_1 =N$.
We shall denote vertical coordinates $\sfx_i \cdot \sfe_2 = u_i$ and 
$\sfy_i\cdot \sfe_2 = v_i$. 
In this way $\eta_i$-s will be modeled  by trajectories $\bbX_i$ of random walks from 
$u_i$ to $v_i$. These walks will be subject to generalized area tilts modulated
by a small parameter $\lambda$. Dyson Ferrari--Spohn diffusions appear in the limit 
$N\to\infty$ and $\lambda\to 0$ subject to the condition in~\eqref{eq:lambdaN-cond} 
below. Let us proceed with precise definitions. 

\smallskip 

\paragraph{\bf Underlying random walk} 

Let $p_y$ be an irreducible random walk kernel on $\bbZ$.
The probability of a finite trajectory $\bbX = (X(-N),X(-N+1),\ldots,X (N))$ is
$\sfp(\bbX) = \prod_i p_{X (i+1) - X (i)}$.
The product probability of $n$ finite trajectories $\underline{\bbX} = 
(\bbX_n,\ldots,\bbX_1)$ is
\be{eq:RW-measure}
\mathbf{P} (\underline{\bbX}) = \prod_{\ell=1}^n \sfp (\bbX_\ell). 
\ee
We shall assume that $\sum_{z\in\bbZ} z p_z = 0$ and that $p$ 
has finite exponential moments; in particular,
\begin{equation}
\label{eq:Assumption}
\sigma^2 \defby  \sum_{z\in\bbZ} z^2 p_z <\infty .
\end{equation}

\paragraph{\bf Generalized area tilts and hard-core constraints}

The hard-core constraint means that we shall consider only ordered 
non-negative trajectories $\underline{\bbX}$. That is, for any $i=-N,\dots,N$, 
the tuple 
\be{eq:hard-core}
\underline{\bbX} (i ) \subset \bbA_n^+\cap \bbZ^n,
\ee
where
\be{eq:Aplus}
\bbA_n^+ = \{ \ur \in \bbR^n \,:\, 0\leq  r_n \leq \dots \leq  r_1\} .
\ee
Given  $\uu,\uv\in\bbZ^n \cap\bbA_n^+$, let $\calP^{\uu,\uv}_{N,+}$ be the family of 
$n$ trajectories $\underline{\bbX}$ starting at $\uu$ at time $-N$, 
ending at $\uv$ at time $N$ and satisfying~\eqref{eq:hard-core}. 

\begin{figure}[t]
\begin{center}
\includegraphics[width=10cm]{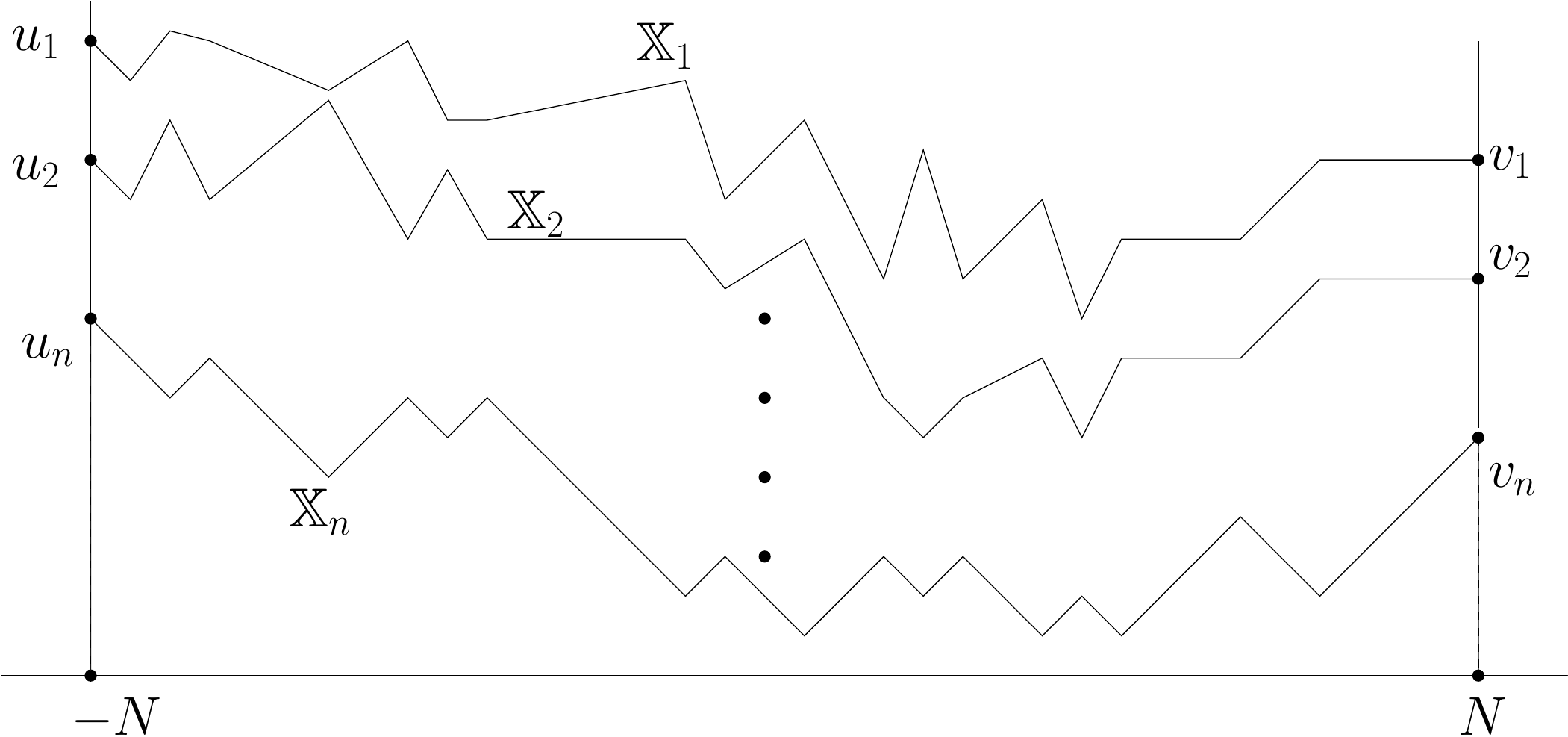}
\end{center}
\caption{The family $\calP^{\uu,\uv}_{N,+}$}
\label{fig:ORWalks}
\end{figure}

Let $\{V_\lambda\}_{\lambda>0}$ be a family of self-potentials 
$V_\lambda: \bbN_0\to\bbR_+$. Given a trajectory $\bbX \subset \bbN_0$ define
\be{eq:Vandw}
 V_\lambda (\bbX )  =\sum_{i=-N}^N V_\lambda ( X (i)) \quad {\rm and} \quad 
 \sfw_\lambda(\bbX) = {\rm e}^{- V_\lambda (\bbX )}\sfp (\bbX ). 
\ee
The term $V_\lambda (\bbX)$ represents a generalized (non-linear)
area below the trajectory $\bbX$. It reduces to (a multiple of) the usual area 
if $V_\lambda(x) = \lambda x$. In general we make the following 
assumptions: 

For any $\lambda>0$, the function $V_\lambda$ on $[0,\infty)$  is 
continuous, monotone increasing and 
satisfies
\be{eq:VL-1}
V_\lambda (0) = 0 \quad\text{ and }\quad \lim_{x\to\infty} V_\lambda(x) = 
\infty .
\ee
In particular, the relation
\be{eq:HL-1}
\Hla^2 V_\lambda(\Hla) = 1
\ee
determines unambiguously the quantity $\Hla$. Furthermore, we make the 
assumptions that $\lim_{\lambda\downarrow 0} \Hla =\infty$ and that there 
exists a function $q\in\sfC^2(\bbR^+)$ such that 
\be{eq:HL-2}
\lim_{\lambda\downarrow 0} \Hla^2 V_\lambda(r\Hla) = q(r),
\ee
uniformly on compact subsets of $\bbR_+$. 
Note  that $\Hla$, respectively $\Hla^2$, plays the role of the spatial, 
respectively temporal, scale in the invariance principle which is formulated 
below in Theorem~\ref{thm:A}. 

A natural class of examples of family of potentials satisfying the above
assumptions is given by $V_\lambda (x ) = \lambda x^\alpha$ with $\alpha>0$.
For the latter, $\Hla=\lambda^{-1/(2+\alpha)}$ and $q(r)=q_0(r)=r^\alpha$.
In this way, the case of linear area tilts $\alpha =1$ corresponds to the
familiar Airy rescaling $H_\lambda = \lambda^{-1/3}$. 

\smallskip 

\paragraph{\bf Probability distributions}

In the case of only one trajectory; $n=1$, given $u,v\in\bbN_0$ and 
$\lambda>0$, define the partition functions and  probability distributions 
\be{eq:pf}
Z^{u,v}_{N,+,\lambda} = 
\sum_{{\bbX}\in\calP^{u,v}_{N,+}} \sfw_\lambda(\bbX)
\text{ and } \bbP^{u,v}_{N,+,\lambda}(\bbX)
= \frac{1}
{Z^{u v}_{N,+,\lambda}}\, \sfw_\lambda(\bbX) 
\1_{\{\bbX\in\calP^{u,v}_{N,+}\} } . 
\ee
In the case of an $n$-tuple $\underline{\bbX} = (\bbX_n,\ldots,\bbX_1)$ of
trajectories, we consider the product weights
$\sfw_\lambda(\underline{\bbX}) = \prod_{i=1}^n \sfw_\lambda(\bbX_i)$
and define the probability distributions $\bbP^{\uu,\uv}_{N,+,\lambda}$
on $\calP^{\uu,\uv}_{N,+}$ by
\be{eq:pd}
\bbP^{\uu,\uv}_{N,+,\lambda} (\underline{\bbX}) 
=
\frac{1}{Z^{\uu,\uv}_{N,+,\lambda}}\,
\sfw_\lambda(\underline{\bbX}) \1_{\{\underline{\bbX} \in 
\calP^{\uu,\uv}_{N,+}\} }  .
\ee

\paragraph{\bf Rescaling and Invariance principle} 

The paths are rescaled as follows: For $t\in H_\lambda^{-2}\bbZ$, define
\be{eq:xN}
\uR^\lambda (t)
= \frac{1}{H_\lambda} \underline{\bbX}  (\Hla^2 t ). 
\ee
Then, extend $\uR^\lambda$ to any $t\in\bbR$ by linear interpolation. 
In this way, given $T>0$ and $\uu, \uv$,
we can talk about the induced distribution
$\bbP^{\uu,\uv;T}_{N,+,\lambda}$
on the space of continuous functions $\sfC\bigl([-T,T],\bbA_n^+\bigr)$.

The  result below  was 
established in \cite{IoShVe2015, IVW16} (under additional mild technical 
assumptions on $V_\lambda$):
\begin{theorem}
\label{thm:A}
For any $n=1, 2, \dots $ fixed, the following happens: 
Let $\lambda_N$ be a sequence satisfying 
\be{eq:lambdaN-cond} 
\lim_{N\to\infty}\lambda_N = 0\quad \text{ and }\quad 
 \lim_{N\to\infty}\frac{N}{H_{\lambda_N}^2} = \infty .
\ee
Fix any $C\in (0,\infty)$ and any $T>0$. 
Then, the sequence of distributions
$\bbP^{\uu,\uv ;T}_{N,+,{\lambda_N}}$ 
converges weakly to the distribution $\bbP_n^{+ ;T}$ of the 
diffusion $\uR (\cdot)$ in~\eqref{eq:DFS-diffusion}, uniformly in $u_1, 
v_1\leq  C H_{\lambda_N}$. 
\end{theorem}

\paragraph{\bf Notes on the proof}

We refer to~\cite{IoShVe2015, IVW16} and to the earlier paper~\cite{HrVe2004} for 
complete statements, further details and careful implementations. Here we shall
try to indicate the main ideas. 

\smallskip 

\noindent
\step{1} Scaling: Consider the probability distribution~\eqref{eq:pf} with 
polymer weights $\sfw_\lambda$ and try to guess what would be the typical 
height $H$ above the wall. There are two factors involved:
\begin{enumerate}[label=(\roman{enumi})]
\item The price for the underlying random walk to stay below height $H$; for this one
pays a constant each $H^2$ units of time.
\item The area tilt which exerts a price $V_\lambda(H)$ each time unit. 
\end{enumerate}
A balance between these two factors is achieved for $H$ satisfying 
$H^2 V_\lambda(H) \sim \mathrm{const}$. Hence the choice~\eqref{eq:HL-1}. 
\smallskip 

\noindent
\step{2} Limiting Sturm--Liouville problem. Consider the partition functions  
$Z^{u,v}_{N,+,\lambda}$ in~\eqref{eq:pf} with the polymer weights
$\sfw_\lambda(\bbX)$ specified in~\eqref{eq:Vandw}. Define 
\be{eq:T-l} 
\sfT_\lambda f (x ) = \sum_{y\in \bbN_0} p_{y-x} {\rm e}^{-V_\lambda (y)} f (y ) .
\ee
In this way,
$Z^{u,v}_{N,+,\lambda} = {\rm e}^{-V_\lambda(u)} (\sfT_\lambda)^{2N}\delta_v (u )$.
Let $f\in \sfC_0^\infty (0, \infty )$. Under the rescaling $x= H_\lambda r$; 
$r\in H_\lambda^{-1}\bbN_0\df \bbN_\lambda$, 
\be{eq:rescT-l} 
\begin{split}
 H_\lambda^2 \lb \sfT_\lambda f (r) - f (r )\rb & = 
 H_\lambda^2 \sum_{s\in \bbN_\lambda} p_{H_\lambda (s-r)}{\rm e}^{-V_\lambda (H_\lambda s)}
 (f (s ) - f (r ))\\ 
 &\hspace{3cm}+ H_\lambda^2\Bigl( \sum_{s\in \bbN_\lambda}
 {\rm e}^{-V_\lambda (H_\lambda s)} - 1\Bigr) f (r )\\
 &\xrightarrow{\eqref{eq:HL-2}, \eqref{eq:Assumption}}
 \frac{\sigma^2}{2}f^{\prime\prime} (r ) - q (r ) f (r ) , 
\end{split}
\ee
which suggests an application of Trotter--Kurtz theorem. 
\smallskip 

\noindent
\step{3} Use of the Karlin--McGregor formula in the case of $n$ ordered walks. 
The Slater determinants~\eqref{eq:DeltaVdm} naturally appear for the limiting ordered 
(continuous) Ferrari--Spohn diffusions via the Karlin--McGregor formula. However, an
application for the discrete partition functions $Z^{\uu,\uv}_{N,+,\lambda}$ in~\eqref{eq:pd}
requires some care. Indeed, the paths can jump over each other, and all the permutation terms
in the general Karlin--McGregor formula~\cite{KMcG} contribute.
One shows, therefore, that as $\lambda\to 0$ the contribution of paths which jump over 
each other vanishes. 

\smallskip 

\noindent
\step{4} Mixing and compactness.
As usual one needs to prove tightness. Furthermore, in order to show that the claim of 
Theorem~\ref{thm:A} holds uniformly in the boundary conditions $\uu$ and $\uv$, one needs
to prove mixing on $\Hla^2$-time scales. These are probabilistic estimates on random walks,
which build upon initial considerations in~\cite{HrVe2004}. 
In the case of $n$ ordered random walks, proving such mixing estimates is by far the most
technically loaded part of~\cite{IVW16}, and it is based on recent strong approximation
techniques developed in~\cite{DW10, DW15, DW15b}.

\appendix

\section{Wulff shapes and Wulff plaquettes}
\label{app:Wulff} 

Let $\varphi_\beta$ be a function on \(\Rd\), which is convex and homogeneous of order one.
We assume that $\varphi_\beta$ is positive on $\Rd\setminus\{0\}$ and that it respects
the $\bbZ^d$-lattice symmetries. 
Depending on the context, the latter may represent either the surface tension
$\varphi_\beta = \tau_\beta$ of low-temperature Ising models, as defined in~\eqref{eq:st-Ising},
or the inverse correlation length $\varphi_\beta = \xi_\beta$ of low-temperature two-dimensional
Ising polymers, as defined in~\eqref{eq:tau-beta}. 

\smallskip 

\paragraph{\bf Wulff shape.} 

By definition $\varphi_\beta$ is the support function of a symmetric compact convex
set $\calW_\beta$ with non empty interior, 
\be{eq:def-Wulff}
\calW_\beta = \bigcap_{\frn\in \bbS^{d-1}}
\setof{x\in\Rd}{x\cdot\frn \leq \varphi_\beta (\frn)}
\ee
In particular, this means that if $\sfn$ is a normal direction to the boundary $\partial  
\calW_\beta$ at $\sfx\in\partial\calW_\beta$, then $\varphi_\beta(\frn) 
= \sfx\cdot\frn$.  

The set $\calW_\beta$ is called the \emph{unnormalized} Wulff shape associated to $\varphi_\beta$.
It is convenient to consider the unit-volume rescaling $\sfW_\beta$ of $\calW$ and the unit-radius
rescaling $\sfK_\beta$ of $\calW$:
\be{eq:rescW} 
\sfW_\beta = \frac{1}{\abs{\calW_\beta}^{1/d}}\calW_\beta\quad {\rm and}\quad 
\sfK_\beta = \frac{1}{\varphi_\beta ({\sfe_1}) } \calW .
\ee
Given a ``sufficiently nice'' (see~\cite{BoIoVe2000} for more precision) subset $V\subset\Rd$,
we define its surface energy
\[
\varphi_\beta({V}) \defby \int_{\partial V} \varphi_\beta(\frn_s) \dd s,
\]
where $\frn_s$ denotes the normal to $\partial V$ at $s$.
We consider the following isoperimetric-type variational problem:

\smallskip
\noindent 
\textbf{(VP)}\quad  Minimize $\varphi_\beta(\partial V)$ among all ``sufficiently 
nice'' sets $V$ of volume $\sfv$,

\smallskip

It can be shown~\cite{Taylor78} that, for any $\sfv>0$, the dilation ${\sfv}^{1/d} \sfW_\beta$
of $\sfW_\beta$ is the unique (up to translation) minimizer of \textbf{(VP)}. 

\smallskip 

In two dimensions, equilibrium crystal shapes associated with surface tension are expected to
have locally analytic and strictly convex boundaries. In the case of the nearest neighbor Ising model,
there is an explicit formula for $\tau_\beta$. However, such a property should follow along the lines
of the Ornstein--Zernike theory whenever phase separation lines admit an effective short-range
random walk representation; see~\cite{CI, CIV08} for such derivations in the context of percolation and
Potts models. 

In higher dimensions, low-temperature Ising Wulff shapes develop facets, at least in lattice
directions~\cite{Miracle99}. Namely, for $\beta$ large, the (sub-differentials) sets
$\partial \tau_\beta (\pm \sfe_i )$ are proper $d-1$-dimensional pieces of $\partial \sfW_\beta$.

\paragraph{\bf Two-dimensional constrained problems and Wulff plaquettes} 

Let us restrict attention to two dimensions. 
In view of the lattice symmetries, $\varphi_\beta(\sfe) \df \varphi_\beta(\pm\sfe_i)$,
where $\pm \sfe_i\; i=1,2$, are unit vectors in lattice coordinate directions, is well defined.
That is, the rescaled shape $\sfK_\beta$ is inscribed into the square $[-1,1]^2$.
Let us define the area $w_\beta$ of $\sfK_\beta$:
\be{eq:W-beta-1} 
w_\beta  = \abs{\sfK_\beta}\in [2,4]. 
\ee
Recall that, given a rectifiable curve $\gamma$, its surface energy is defined by 
$\varphi_\beta(\gamma) = \int_\gamma \varphi_\beta(\frn_s)\dd s$. 
Consider the following \emph{constrained} isoperimetric problem~\cite{SS96}:

\smallskip
\noindent 
\textbf{(CVP)}\quad  Minimize $\varphi_\beta( V)$ among all subsets $V\subset [-1, 1]^2$
with rectifiable boundary $\partial V$ and area $\abs{V} =\sfb \in (0,4]$. 
\smallskip

The answer depends on the value of $\sfb$. It is natural to state it in terms of the unit-radius
shape $\sfK_\beta$. 

Since $\varphi_\beta(\cdot)$ is the support function of $\calW_\beta$, the
ratio $\frac{\varphi_\beta(\cdot)}{\varphi_\beta(\sfe)}$ is the support function of $\sfK_\beta$. Hence, 
\be{eq:tension-Wb} 
w_\beta 
=
\int_{\sfK_\beta}\dd \sfx
=
\frac{1}{2}\int_{\sfK_\beta} \mathrm{div} (\sfx) \dd\sfx
=
\frac{1}{2} 
\frac{\varphi_\beta (\partial\sfK_\beta)}{\varphi_\beta (\sfe )}
\quad \Rightarrow\quad
\varphi_\beta (\sfK_\beta)
=
2\varphi_\beta (\sfe ) w_\beta .
\ee
In the above notation, the Wulff shape of area $\sfb$ is
$\sfK_\beta^\sfb = \sqrt{\frac{\sfb}{w_\beta}}\sfK_\beta$
and, for $\sfb\in (0, w_\beta]$,   
\be{eq:Wulff-b} 
\varphi_\beta (\sfb ) \df 
\mintwo{V~:~ \abs{V} = \sfb}{V\subseteq [-1,1]^2}
\varphi_\beta( V )
=
\varphi_\beta( \sfK_\beta^\sfb)
= 
2\varphi_\beta (\sfe ) \sqrt{{\sfb}{w_\beta}}. 
\ee
Clearly, $\sfK_\beta^\sfb$ is the unique (again up to translations within $[-1,1]^2$)
solution of \textbf{(CVP)} whenever $\sfb \in (0,w_\beta]$. Note that the radius $r_{\sfb}$ of $\sfK_\beta^\sfb$
satisfies 
\be{eq:rb-K} 
 r_{\sfb} = \sqrt{\frac{\sfb}{w_\beta}}\quad {\rm and}\quad 
 \frac{\dd \varphi_\beta (\sfK_\beta^\sfb )}{\dd \sfb}  = 
\frac{\varphi_\beta (\sfe )}{r_\sfb } .
\ee
For $\sfb \in (w_\beta , 4]$, the shape $\sfK_\beta^\sfb$ does not fit into the square $[-1, 1]^2$.
In this case~\cite{SS96},
\be{eq:Wulff-P-b} 
\varphi_\beta (\sfb )\df \mintwo{V~:~ \abs{V} = \sfb}{V\subseteq [-1,1]^2}
\varphi_\beta(V )
= 
\varphi_\beta( \sfP^\sfb_\beta)
=
\varphi_\beta (\sfe )\cdot \Bigl( 8 - 2\sqrt{(4-w_\beta)(4-\sfb)} \Bigr)  .
\ee
Above, the \emph{Wulff plaquettes} $\sfP^\sfb_\beta$ are constructed as follows: Place four Wulff
shapes of radius 
\be{eq:rb}
r_\sfb =\sqrt{\frac{4-\sfb}{4-w_\beta}}\leq 1 
\ee
into the four corners of $[-1,1]^2$, in such a way that each of these shapes is tangent to the two
corresponding sides of the square.
Then take the convex envelope. Note that, by construction, $\partial\sfP^\sfb_\beta$ has four flat
segments of length $2(1-r_\sfb)$ on each of the four sides of $[-1,1]^2$. Also, note that in terms
of $r_\sfb$. The area $\sfb$ and the surface energy $\varphi_\beta (\sfP_\beta^\sfb)$ of the Wulff
plaquette $\sfP_\beta^\sfb$ can be recorded as
\be{eq:rb-b-Pb} 
\sfb = 4 - (4-w_\beta )r_\sfb^2 ,\ 
\varphi_\beta ( \sfP_\beta^\sfb ) = \varphi_\beta (\sfe )\lb 
8 -  2 (4 -w_\beta )r_\sfb\rb \  {\rm and}\ 
\frac{\dd \varphi_\beta (\sfP_\beta^\sfb )}{\dd \sfb}  = 
\frac{\varphi_\beta (\sfe )}{r_\sfb }, 
\ee
where the latter identity directly follows from~\eqref{eq:Wulff-P-b} and~\eqref{eq:rb}. 

\smallskip 

\paragraph{\bf Dual constrained variational problems}
Let $\nu>0$ and consider:  

\smallskip
\noindent 
\textbf{(DCVP)}\quad 
Find\qquad 
$\max_{\sfa\in [0,4]} \{\nu \sfa - \varphi_\beta(\sfa)\}$. 

\smallskip 

By~\eqref{eq:Wulff-b} and~\eqref{eq:Wulff-P-b}, the function $\varphi_\beta$ is 
concave on $[0, w_\beta]$ and convex on $[w_\beta,4]$. Therefore, there 
exists a critical value $\nu_\beta \in (0, \infty)$ such that $\sfa=0$ is 
the unique solution to~\textbf{(DCVP)} for $\nu\in [0, \nu_\beta )$, while there is a
unique solution $\sfa = \sfa(\beta,\nu)\in (w_\beta,4)$ for every $\nu>\nu_\beta$.
At $\nu = \nu_\beta$, there are two solutions: $\sfa=0$ and
$\sfa = \sfa_\beta > w_\beta$. The critical pair $(\nu_\beta,\sfa_\beta)$ satisfies
the following equation: 
\be{eq:eq-for-nu} 
\varphi_\beta^{\prime} (\sfa_\beta )
=
\varphi_\beta (\sfe )\cdot \sqrt{\frac{4-w_\beta}{4-\sfa_\beta}}
=
\nu_\beta
\text{ and }
\nu_\beta \sfa_\beta
= 
\varphi_\beta (\sfe )\cdot  \Bigl( 8 -  2\sqrt{(4-w_\beta)(4-\sfa_\beta)} \Bigr) .
\ee
Geometrically, if one starts to rotate counter-clockwise a line passing through zero,
then the latter will touch for the first time the graph of
$\sfa \mapsto \varphi_\beta (\sfa)$ at the angle $\nu_\beta$ and at the point
$(\sfa_\beta,\varphi_\beta(\sfa_\beta))$. 
 
Note that, for $\nu >\nu_\beta$, the optimal shape is necessarily a Wulff 
plaquette and that there is a strict inclusion of optimal shapes that 
correspond to different $\nu^\prime > \nu >\nu_\beta$.  

In order to get an explicit expression for the critical slope $\nu_\beta$,
it is convenient to compare~\eqref{eq:eq-for-nu} with~\eqref{eq:rb-b-Pb}. 
Let $\rho_\beta$ be the critical radius. Then the first, and consequently 
the second,  of \eqref{eq:eq-for-nu} read as 
\be{eq:eq-for-nu-r} 
\nu_\beta = \frac{\varphi_\beta (\sfe )}{\rho_\beta}\quad {\rm and}\quad 
\frac{4 - (4-w_\beta )\rho_\beta^2}{\rho_\beta} = 8 - 2 (4 -w_\beta )\rho_\beta .
\ee 
Solving the second of~\eqref{eq:eq-for-nu-r} (and taking into account that
the solution we need should satisfy $\rho_\beta \leq 1$), we infer: 
\be{eq:crit-beta} 
\rho_\beta = \frac{2}{2 +\sqrt{w_\beta}}\ \Rightarrow\ 
\nu_\beta = \frac{\varphi_\beta (\sfe )\lb 2 + \sqrt{w_\beta}\rb}{2} = 
\frac{4\varphi_\beta (\sfe ) + \varphi_\beta (\sfW_\beta )}{4} , 
\ee
where we relied on~\eqref{eq:Wulff-b} to identify the energy of the
unit-volume Wulff shape as
$\varphi_\beta (\sfW_\beta )=2\varphi_\beta (\sfe )\sqrt{w_\beta}$.

\section{Ferrari--Spohn diffusions}
\label{app:FS}

In this section, the notations $\|\cdot\|_2$ and $\langle \cdot,\cdot\rangle_2$ are
reserved for the norm and scalar product in $\bbL_2(\bbR_+)$.

Given  $\sigma>0$ and  a non-negative function $q\in\sfC^2(\bbR_+)$ which 
satisfies $\lim_{r\to\infty} q (r) = \infty$, consider the following family 
of singular  Sturm--Liouville operators on
$\bbR_+$:
\begin{equation}
\label{eq:SL-operators}
\sfL_{\sigma,q} = \frac{\sigma^2}{2}\frac{\dd^2}{\dd r^2} - q (r) .
\end{equation}
with boundary condition $\varphi(0) = 0$. $\sfL_{\sigma,q}$ possesses a 
complete orthonormal family $\lbr \varphi_i\rbr$  of simple eigenfunctions
in $\bbL_2(\bbR_+)$ with eigenvalues
\begin{equation}
\label{eq:eigenv-ST}
0 >- \eig_0 > -\eig_1 > -\eig_2 > \dots ;\ \lim \eig_j = \infty.
\end{equation}
The eigenfunctions $\varphi_i$ are smooth and $\varphi_i$ has exactly $i$
zeros in $(0,\infty)$, $i=0,1,\ldots $.

The Ferrari--Spohn diffusion, associated to \(\sigma\) and \(q\), is the 
diffusion on $(0,\infty)$ with generator
\be{eq:FS-Gen}
\sfG_{\sigma, q}\psi
\df
\frac{1}{\varphi_0 } (\sfL_{\sigma,q} +\eig_0) (\psi\varphi_0) =
\frac{\sigma^2}{2}\frac{\dd^2\psi }{\dd r^2} + \sigma^2
\frac{\varphi_0^\prime }{\varphi_0}\frac{\dd\psi }{\dd r} =
\frac{\sigma^2}{2\varphi_0^2}
\frac{\dd }{\dd r}\Bigl(\varphi_0^2\frac{\dd \psi }{\dd r}\Bigr).
\ee
This diffusion is ergodic and reversible with respect to the measure
$\dd\mu_0(r)=\varphi_0^2(r)\dd r$. We denote by $\sfS_{\sigma, q}^t$ the 
corresponding semigroup and by $\bbP_{\sigma,q}$ the associated path 
measure.

\medskip
The most relevant case for the present paper is when the function $q$ is 
linear, \(q(r)= cr\). In that case, since the Airy function \(\sf Ai\) 
satisfies \(\frac{\dd^2}{\dd r^2}{\sf Ai} (r) = r{\sf Ai}(r)\), one can 
easily check that
\be{eq:AiryGen}
\varphi_0 = {\sf Ai}(\chi r-\omega_1) \quad{\rm and}\quad e_0 =
\frac{c\omega_1}{\chi},
\ee
where $-\omega_1 = -2.33811\ldots$ is the first zero of ${\sf Ai}$ and
$\chi = \sqrt[3]{\frac{2c}{\sigma^2}}$.

\smallskip 

\paragraph{\bf Dyson Ferrari--Spohn diffusions} 

Let us fix  $n\in\bbN$. Recall  the notation 
\[
\bbA_n^+ = \{ \ur \in \bbR^n \,:\, 0 \leq  r_n \leq \dots \leq  r_1\} .
\]
Let $\uX(t) = (\sfX_n (\cdot), \dots, \sfX_1(\cdot))$ be $n$ independent
copies of Ferrari--Spohn diffusions starting at $\ur \in \bbA_n^+$. 
Define the killing time  
\be{eq:Tkilled}
\tau  = \min\{ t \,:\, \uX(t)\not\in\bbA_n^+ \} .
\ee
In other words, $\tau$ is the minimum between the first collision time and
the first time the bottom trajectory exits from the positive semi-axis.
Dyson Ferrari--Spohn diffusion $\uR(t)$ is $\uX(\cdot)$ conditioned to 
survive forever, that is, conditioned on $\{\tau = \infty\}$. We refer
to~\cite{IVW16} for a precise construction. As is explained in the latter
paper, a Dyson Ferrari--Spohn  diffusion $\uR(t)$ satisfies the following SDE: 
Let
\be{eq:DeltaVdm}
\Delta (\ur) 
= 
\operatorname{det} 
\begin{bmatrix}
\varphi_1(r_1)	& \varphi_2(r_1) 	& \cdots &\varphi_n (r_1 ) \\
\varphi_1(r_2)	& \varphi_2(r_2) 	& \cdots &\varphi_n (r_2 ) \\
\vdots 		& \vdots 		& \ddots &\vdots \\
\varphi_1(r_n)	& \varphi_2(r_n) 	& \cdots &\varphi_n (r_n ) 
\end{bmatrix} .
\ee
be the Slater determinant of $\sfL_{\sigma,q}$. Alternatively, $\Delta(\ur)$
is the first Dirichlet eigenfunction of $\sfL_n +\dots + \sfL_1$
on $\bbA_n^+$, where $\sfL_i$ is a copy of $\sfL_{\sigma,q}$ acting on the $i$-th
coordinate $r_i$\footnote{We refer to~\cite{Sosh00} for an exposition of
fermionic ground states and determinantal point processes}. Then,  
\be{eq:DFS-diffusion}
\dd\uR (t) = \sigma \dd\uB (t) + \nabla\log(\Delta) (\uR(t)) \,\dd t ,
\ee
where $\uB(t)$ is the standard $n$-dimensional Brownian motion. 
Note that the diffusion $\uR(\cdot)$ on $\bbA_n^+$ is ergodic and reversible
with respect to $\Delta^2(\ur)\,\dd\ur$.
We use $\bbP_n^{+ ;T}$ to denote its (stationary) distribution on the space of
continuous functions $\sfC\bigl([-T,T],\bbA_n^+\bigr)$. 

\section{Self-avoiding walks under area tilts} 
\label{app:SAW} 

In this appendix, we shall sketch how Ferrari--Spohn diffusions appear 
as scaling limits of super-critical nearest-neighbor 
two-dimensional (vertex) self-avoiding 
walks (SAW-s) under area tilts. Being an instance of Ising polymers without 
cluster decorations super-critical SAW-s are viewed as simplified 
models of phase separation lines in 2D Ising models below critical 
temperature, and the discussion below is intended to reinforce 
Conjecture~\ref{conj:IsingFS}. A complete proof will hopefully appear elsewhere. 

Let us fix $\beta > \log \mu_c$, where ${\mu_c}$ is the connectivity constant of 
$\bbZ^2$. In the sequel, we shall drop $\beta$ from all the notation. The reference weight
of SAW trajectory $\gamma = \lb \gamma (0 ), \dots , \gamma (n )\rb$ is 
${\rm e}^{-\beta n} = {\rm e}^{-\beta \abs{\gamma}}$. A path $\gamma$ is said to be a 
positive bridge, $\gamma\in\calB^+$, if $\gamma \subset \bbH_+^2$ and 
if $\gamma (0)\cdot \sfe_1 \leq \gamma (\ell )\cdot \sfe_1 \leq \gamma (n )\cdot \sfe_1 $.
Given $x, y\in \bbN$, define 
\be{eq:BplusN} 
\calB^+_N (x, y ) = \setof{ \gamma \in \calB^+}{\gamma (0 ) = (-N , x)\ {\rm and}\ \gamma (n ) = (N , y)} .
\ee
Note that $\abs{\gamma} = n$ is not fixed in the above definition; we just talk about all the
bridges from $(-N , x)$ to $(N , y)$. For any $\gamma\in \calB^+$, the area $A (\gamma )$ 
is well defined. Given $\lambda >0$ (area tilt), consider the weights $w_\lambda$ and the
probability distributions $\bbP_{N, +, \lambda}^{x, y}$: 
\be{eq:wl-Pplus} 
w_\lambda (\gamma ) = {\rm e}^{-\beta \abs{\gamma} - \lambda A(\gamma )}\quad {\rm and}\quad
\bbP_{N, +, \lambda}^{x, y} (\gamma ) = \frac{w_\lambda (\gamma )}{Z_{N, +, \lambda}^{x, y}}
\1_{\lbr \gamma\in \calB^+_N (x, y )\rbr}.
\ee
We shall try to explore what happens when $\lambda\to 0$ and $N\to\infty$ in an appropriate
way. 
The inverse correlation length $\xi_\beta$ is defined and positive for any 
$\beta > \log\mu_c$, and Ornstein--Zernike theory, as described in Subsection~\ref{sub:Ising}, 
applies \cite{CC86, I98, IV08}.  Hence,
for $\lambda $  sufficiently small and $N$ sufficiently large,
we can restrict attention to 
$\gamma$ which admit an irreducible 
decomposition \eqref{eq:a-concat} with respect to some symmetric cone 
$\calY_\sfh$ along 
the horizontal axis, where $\sfh = \sfh_\beta = \nabla \xi_\beta (\sfe_1 ) = \xi_\beta (\sfe_1 ) \sfe_1$. 
If $\gamma$ admits \eqref{eq:a-concat}, then we can think about linear interpolation 
$\hat{\gamma}$ through the vertices of \eqref{eq:a-concat}. Notice that 
$\hat{\gamma}$ is already a well-defined function on $[-N , N]$. Let us rescale it 
as 
\be{eq:x-lambda} 
x_\lambda (t ) = \lambda^{1/3} \hat{\gamma} (\lambda^{-2/3}t ).
\ee
Under $\bbP_{N, +, \lambda}^{x, y}$, the rescaled path   $x_\lambda (\cdot )$ is 
viewed as a random continuous function on $[-N \lambda^{2/3},N \lambda^{2/3}]$. 
Here is the SAW counterpart of Conjecture~\ref{conj:IsingFS}: 
\begin{theorem} 
 \label{thm:SAW-FS} 
Consider $\lambda = \lambda_N$ which satisfies $\lim\lambda_N = 0$ and 
$\lim N\lambda_N^{2/3} = \infty$. Consider a sequence $\lbr x_N , y_N\rbr$ 
of boundary conditions, such that $x_N, y_N \leq C\lambda_N^{-1/3}$. 

Then the distribution of $x_{\lambda_N} (\cdot )$ 
under $\bbP_{N, +, \lambda}^{x_N, y_N}$ is weakly convergent to the distribution of 
FS diffusion with parameters $\sigma^2 = \chi_\beta$ and $q (r ) = r$, where 
$\chi_\beta$ is the curvature of $\partial\calW_\beta$ at $\sfh_\beta$. 
\end{theorem}
Below we indicate main steps of the proof. 
\smallskip 

\noindent 
\step{1} Identification of $\chi_\beta$ in terms of the probability distribution 
$\bbP_\beta^\sfh$ (compare with  \eqref{eq:OZ-weights}) on the set of irreducible 
animals/paths $\sfD \df \sfD_\sfh$: 
\be{eq:OZ-weight-SAW} 
 \bbP_\beta^\sfh (\fra ) = {\rm e}^{\sfh\cdot \sfX (\fra ) - \beta\abs{\gamma }} .
\ee
By \eqref{eq:p-Wbeta-eq} the boundary $\partial\calW_\beta$ is parametrized around 
$\sfh$ as $t\to \phi (t )$ with 
\be{eq:param-Wb} 
 \bbE_\beta^\sfh \lb {\rm e}^{ (-\phi (t ) , t )\cdot \sfX }\rb   \equiv 1 .
\ee
Writing $\sfX = (\sfT , \sfY )$ one, using symmetry of $\sfY$ under $\bbP_\beta^\sfh$ 
and 
second order expansion, recovers $\chi_\beta$ as
\be{eq:chi-b-var} 
\chi_\beta = \frac{\bbE_\beta^\sfh (\sfY^2 )}{\bbE_\beta^\sfh (\sfT )} .
\ee
\step{2} Reduction to effective Markov chains. Instead of $\calB_N^+$, consider the
concatenations of, say, $m$ irreducible pieces. It happens to be convenient to 
record such concatenation as 
\be{eq:cat-S} 
\lbr \fra_1 , x_1\rbr\circ \lbr \fra_2 , x_2\rbr\circ\cdots\circ  \lbr \fra_m , x_m\rbr .
\ee 
We understand a pair $\{\fra,x\}$ as an animal $\fra$ being attached to a point 
at height $x\in\bbN$ from the left,  
and we write $\{\fra,x\} \in \calS_+$ if 
the resulting structure is inside $\bbH_+^2$. Let us record  the area 
\be{eq:area-pair}
 A\lb \lbr \fra , x\rbr \rb = \lb x - \sfY (\fra )\rb \sfT (\sfa ) + \Delta (\fra ) , 
\ee
where it should be understood that the above relation defines $\Delta (\fra )$. 

\begin{center}
\label{fig:animal-SAW}
\includegraphics[width=3.5cm]{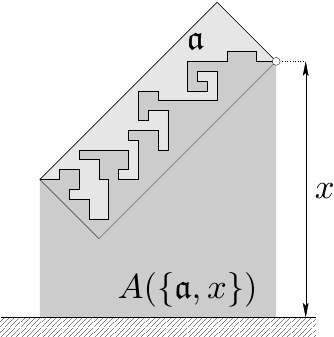}\\
\centerline{A pair $\{\fra,x\} \in \calS_+$ and its area $A (\{\fra,x\})$}
\end{center}

Finally define a sub-Markov kernel $\sfQ_\lambda$ on $\calS_+$: 
\be{eq:Q-lambda} 
Q_\lambda\lb \lbr \fra , u\rbr , \lbr \frb , v\rbr\rb = \1_{\lbr u = v - \sfY (\frb )\rbr}
\bbP_\beta^{\sfh} (\frb ) {\rm e}^{- \lambda A\lb \lbr \frb , v\rbr \rb } .
\ee
In this way, we can record various partition functions
in terms of powers of $Q_\lambda$.
For instance,
the partition function of all $m$-irreducible step trajectories from height $x$ 
to height $y$  reads as:
 \be{eq:pf-powers}
 \sumtwo{\lbr \fra , u\rbr \in\calS_+}{u - \sfY (\fra ) =x}
 \sum_{\frb :\lbr \frb , y \rbr \in\calS_+ } 
 Q_\lambda^m\lb \lbr \fra , u\rbr , \lbr \frb , y \rbr\rb
\ee

\begin{center}
\includegraphics[width=\textwidth]{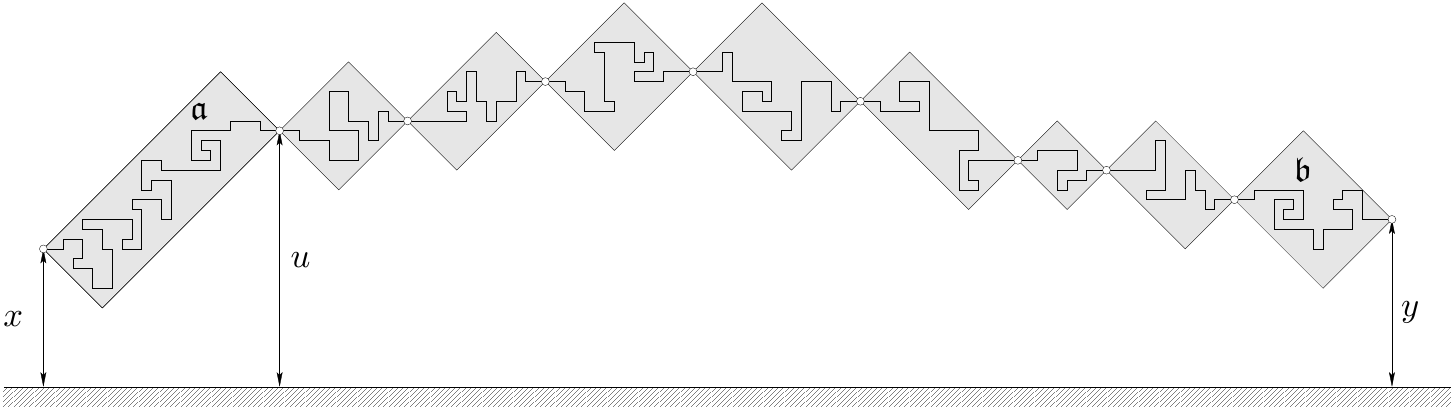}\\
\centerline{Path of irreducible steps from $x$ to $y$.}
\end{center}

\step{3} Rescaling and application of Trotter--Kurtz theorem. 
As before, define $\bbN_\lambda = \lambda^{1/3}\bbN$. For $r, s\in \bbN$, consider the 
rescaled sub-Markov kernel
\be{eq:hatQlambda}
\hat{Q}_\lambda (\{\fra,r\}, \{\frb,s\}) = 
Q_\lambda (\{\fra,r\lambda^{-1/3}\}, \{\frb,s\lambda^{-1/3}\}).
\ee
Let $f\in \sfC_0^\infty \lb 0, \infty\rb$. Given $r\in\bbN_\lambda$ and $\fra$ such that
$\lbr\fra, \lambda^{-1/3} r\rbr\in \calS_+$, notice that
\begin{align}
\label{eq:Ql-gen}
\frac{1}{\lambda^{2/3}}& \Bigl( \sum_{\lbr \frb , s\rbr} 
\hat{Q}_\lambda \lb \lbr \fra , r\rbr , \lbr \frb , s \rbr\rb f (s ) - f (r )\Bigr)\\
&\stackrel{\eqref{eq:Q-lambda}}{=} \frac{1}{\lambda^{2/3}} \Bigl( \sum_{\frb\in \sfD}
\1_{\lbr \lambda^{-1/3}r +\sfb \subset \bbH_+^2\rbr} \bbP_\beta^h (\frb ) 
{\rm e}^{-\lambda (\lambda^{-1/3} r \sfT (\frb ) +\Delta (\frb ))} f (r + \lambda^{1/3}
\sfY (\frb )) - f (r )\Bigr) \\ 
& =\frac{\bbE_\beta^\sfh (\sfY^2 ) }{2} f^{\prime\prime} (r ) - \bbE_\beta^\sfh (\sfT ) r
f (r ) +\smo{1} \stackrel{\eqref{eq:chi-b-var}}{=} 
\bbE_\beta^\sfh (\sfT )\bigl( \frac{\chi_\beta}{2}f^{\prime\prime} (r ) - r f (r ) 
+\smo{1} \bigr)\\ 
& \df \bbE_\beta^\sfh (\sfT )\sfL_\beta f +\smo{1} .
\end{align}
This, when relying on Theorem~I.6.5 in~\cite{EK}\footnote{See Section~2.2 in~\cite{IoShVe2015} 
for more details in the case usual random walks}, means that, for any $t >0$, any 
$g, f\in \sfC_0^\infty (0, \infty )$ and any $\fra \in \sfD$,
\be{eq:TK-lim} 
\lim_{\lambda\to 0} \lambda^{1/3} \sum_{r, s\in \bbN_\lambda}\sum_{\frb\in \sfD} 
g (r ) \hat{Q}_\lambda^{\lfloor t \lambda^{-2/3}/\bbE_\beta^\sfh (\sfT )\rfloor} 
\lb\lbr \fra , r\rbr , \lbr \frb , s \rbr\rb f (s) = 
\int_0^\infty f (r ){\rm e}^{t \sfL_\beta } g(r )\dd r .
\ee
\step{4} Mixing and tightness for the rescaling~\eqref{eq:x-lambda} 
of $\hat{\gamma}$ follows by arguments similar to those employed in the 
case of the usual random walk in~\cite{IoShVe2015}. 
Mixing is needed for coupling distributions with 
fixed horizontal projections as in \eqref{eq:wl-Pplus} with distributions which 
correspond to a fixed number of irreducible steps, as generated by the sub-Markov
kernel~\eqref{eq:Q-lambda}, as well as for coupling distributions with different numbers
of irreducible steps and with the limiting ergodic Markov chain on $\sfS_+$, which 
is defined for any $\lambda >0$. Once mixing is established, \eqref{eq:TK-lim} implies
convergence of the finite-dimensional distributions of $x_\lambda $ to the finite-dimensional
distributions of the FS $(\chi_\beta , r)$-diffusion, which is upgraded to a full invariance
principle by  tightness.

\bibliographystyle{plain}
\bibliography{MPRF.bib}

\end{document}